\documentclass[preprint,12pt]{elsarticle}
\usepackage{amssymb}

\journal{Journal of Approximation Theory}

\usepackage{graphicx,amsmath,amsfonts}
\usepackage{verbatim}  
\usepackage{pict2e}

\usepackage{times,amsfonts,amssymb,amsmath}
\usepackage{psfrag, color}
\usepackage{multirow}
\usepackage{float}
\usepackage{enumitem}
\def\NN{\mathbb{N}}
\def\RR{\mathbb{R}}

\def\sign{{\rm sign}}

\newcommand{\etc}{\emph{etc}}
\DeclareMathOperator\supp{supp}

\def\bfms#1{{\boldsymbol{#1}}}
\newcommand{\pika}[0]{{\raise0.5ex\hbox{{\bf .}}}}
\newcommand{\mybound}[2]{\genfrac{}{}{0pt}{}{#1}{#2}}
\newcommand{\etal}{\emph{et~al.}}
\newcommand{\figGeomIntLeft}{figGeomIntLeft.eps}
\newcommand{\figGeomIntRight}{figGeomIntRight.eps}
\newcommand{\figNorDiffLeft}{figNorDiffLeft.eps}
\newcommand{\figNorDiffRight}{figNorDiffRight.eps}
\newcommand{\figSdefL}{figSdefL.eps}
\newcommand{\figLemcKKpOnePf}{figLemcKKpOnePf.eps}
\newcommand{\figgBPMaL}{figgBPMaL.eps}
\newcommand{\figgBPMaR}{figgBPMaR.eps}
\newcommand{\figgBPMbL}{figgBPMbL.eps}
\newcommand{\figgBPMbR}{figgBPMbR.eps}
\newcommand{\figLemKKpiPfL}{figLemKKpiPfL.eps}
\newcommand{\figLemKKpiPfR}{figLemKKpiPfR.eps}
\newcommand{\figLemabKKpOnePfL}{figLemabKKpOnePfL.eps}
\newcommand{\figLemabKKpOnePfR}{figLemabKKpOnePfR.eps}

\newcommand{\figRemNisTwoL}{figRemNisTwoL.eps}
\newcommand{\figRemNisTwoR}{figRemNisTwoR.eps}
\newcommand{\figApprOrdLeft}{figApprOrdLeft.eps}
\newcommand{\figApprOrdRight}{figApprOrdRight.eps}
\newcommand{\figLemabKKpiPfL}{figLemabKKpiPfL.eps}
\newcommand{\figLemabKKpiPfR}{figLemabKKpiPfR.eps}
\newcommand{\fignIsthreeSlice}{fignIs3Slice.eps}
\newcommand{\fignIsthreeL}{fignIs3L.eps}
\newcommand{\figSjmOnejjPOne}{figSjmOnejjPOne.eps}
\newcommand{\figSjmOnejjPTwo}{figSjmOnejjPTwo.eps}
\newcommand{\refref}[2]{(#1~\ref{#2})}
\newcommand{\refnref}[2]{#1~\ref{#2}}
\newcommand{\figref}[1]{\refref{Fig.}{#1}}
\newcommand{\fignref}[1]{\refnref{Fig.}{#1}}

\newcommand{\lemnref}[1]{\refnref{Lemma}{#1}}

\newcommand{\thmnref}[1]{\refnref{Theorem}{#1}}
 
\newcommand{\remnref}[1]{\refnref{Remark}{#1}}

\newcommand{\ennref}[1]{\textit{\ref{#1}}}

\newcommand{\interior}[1]{\textrm{int}\left( #1 \right)}
\newcommand{\vv}[1]{\bfm{#1}}
\newcommand{\lin}{ \mathcal{L}_{\vv{t},\vv{c},j}}
\newcommand{\linp}[1]{ \lin\left(#1\right) }
\newcommand{\vdva}[2]{ \begin{pmatrix} #1 \\ #2 \end{pmatrix}}
\newcommand{\boud}{\partial}
\newcommand{\boudp}[1]{\boud\,#1}
\newcommand{\vcc}[3]{\left( #1 \right)_{#2}^{#3}}
\newcommand{\partName}{t}
\newcommand{\vt}{\bfms{\partName}}

\newcommand{\vtZ}{\vt^0}
\newcommand{\vtj}[1]{\vt_{\backslash\left\{#1\right\}}}

\newcommand{\vtZj}[1]{\vtZ_{\backslash\left\{#1\right\}}}
\newcommand{\ti}[1]{\partName_{#1}}
\newcommand{\tti}[1]{\widetilde{\partName}_{#1}}
\newcommand{\tiZ}[1]{\partName^0_{#1}}

\newcommand{\triang}{\ensuremath{\mathcal{T}}}
\newcommand{\qria}{\bfm{V}}
\newcommand{\qriac}[1]{\qria_{#1}}
\newcommand{\qrial}{\qriac{0}}
\newcommand{\qriar}{\qriac{3}}
\newcommand{\qriaur}{\qriac{2}}
\newcommand{\qriaul}{\qriac{1}}

\newcommand{\tria}{\bfm{T}}
\newcommand{\triac}[1]{\tria_{#1}}
\newcommand{\trial}{\triac{0}}

\newcommand{\triaur}{\triac{2}}
\newcommand{\triaul}{\triac{1}}
\newcommand{\charfunName}{\chi}
\newcommand{\charfun}{\charfunName}
\newcommand{\charfuni}[1]{\charfun_{#1}}
\newcommand{\charfunParName}{\eta}
\newcommand{\charfunPara}{\charfunParName}

\newcommand{\charfunParb}{\widetilde{\eta}}

\newcommand{\charfunNopp}[2]{\charfun_{{_{#1}}}\left(#2\right)}
\newcommand{\charfunp}[1]{\charfun_{{_{#1,\charfunPara,\charfunParb}}}}
\newcommand{\charfunpp}[2]{\charfunp{#1}\left(#2\right)}
\newcommand{\tht}[2]{\theta_{#1,{\scriptscriptstyle#2}}}
\newcommand{\thta}[1]{\tht{\vv{a}}{#1}}
\newcommand{\thtb}[1]{\tht{\vv{b}}{#1}}
\newcommand{\thtc}[1]{\tht{\vv{c}}{#1}}
 \newcommand{\thtanic}{\thta{\trial}}
 \newcommand{\thtadva}{\thta{\triaur}}
 \newcommand{\thtbnic}{\thtb{\trial}}
 \newcommand{\thtbdva}{\thtb{\triaur}}
 \newcommand{\thtcnic}{\thtc{\trial}}
 \newcommand{\thtcdva}{\thtc{\triaur}}
\newcommand{\ctau}[2]{\tau_{#1,{\scriptscriptstyle#2}}}
\newcommand{\taua}[1]{\ctau{\vv{a}}{#1}}
\newcommand{\taub}[1]{\ctau{\vv{b}}{#1}}
\newcommand{\tauc}[1]{\ctau{\vv{c}}{#1}}
 \newcommand{\tauanic}{\taua{\qrial}}
 \newcommand{\tauaena}{\taua{\qriaul}}
 \newcommand{\tauadva}{\taua{\qriaur}}
 \newcommand{\tauatri}{\taua{\qriar}}
 \newcommand{\taubnic}{\taub{\qrial}}
 \newcommand{\taubena}{\taub{\qriaul}}
\newcommand{\taubdva}{\taub{\qriaur}}
\newcommand{\taubtri}{\taub{\qriar}}

 \newcommand{\taucena}{\tauc{\qriaul}}

\newcommand{\tauctri}{\tauc{\qriar}}

\newcommand{\sepsname}{\mathcal{S}}
\newcommand{\sm}{\sepsname}

\newcommand{\co}{\textrm{co}}
\newcommand{\cop}[1]{\co\!\left\{#1\right\}}
\newcommand{\smfp}[1]{\sm_{\backslash\left\{#1\right\}}}
\newcommand{\smint}{\interior\sm}

\newcommand{\boundary}[1]{\partial #1}

\newcommand{\deter}[3]{\begin{Vmatrix} #1 \end{Vmatrix}_{#2}^{#3} }
\newcommand{\da}{\Delta{a}}
\newcommand{\deap}[1]{\da_{#1}}
\newcommand{\db}{\Delta{b}}
\newcommand{\debp}[1]{\db_{#1}}
\newcommand{\dc}{\Delta{c}}
\newcommand{\decp}[1]{\dc_{#1}}
\newcommand{\dta}{\Delta\widetilde{a}}
\newcommand{\detap}[1]{\dta_{#1}}
\newcommand{\dtb}{\Delta{\widetilde{b}}}
\newcommand{\detbp}[1]{\dtb_{#1}}

\newcommand{\om}[4]{\prod\limits_{\mybound{\ell = #1}{\ell \ne#3}}^{#2}\left(#4 - t_{\ell}\right)}

\newcommand{\omj}[4]{\prod\limits_{\mybound{\ell = #1}{\ell \ne#3}}^{#2}\left(#4 - t_{j_\ell}\right)}

\newcommand{\omd}[1]{\omega_{#1}}
 \newcommand{\omdp}[2]{\omd{#1}\left( #2 \right)}
\newcommand{\omdcp}[2]{\omd{#1}'\left( #2 \right)}
\newcommand{\fun}{f}

\newcommand{\funpp}[2]{\fun_{#1,#2}}

\newcommand{\fa}{\funpp{\bfm{a}}{j}}
\newcommand{\faap}[1]{\funpp{\bfm{a}}{#1}}
\newcommand{\fb}{\funpp{\bfm{b}}{j}}
\newcommand{\fbbp}[1]{\funpp{\bfm{b}}{#1}}
\newcommand{\fc}{\funpp{\bfm{c}}{j}}
\newcommand{\fap}[1]{\fa\left( #1 \right)}
\newcommand{\fbp}[1]{\fb\left( #1 \right)}
\newcommand{\fcp}[1]{\fc\left( #1 \right)}

\newcommand{\var}{\mathcal{V}}
\newcommand{\vari}[1]{\ensuremath{\var_{#1}}}
\newcommand{\varip}[2]{\ensuremath{\vari{#1}\left(#2\right)}}
\newcommand{\varisp}[1]{\varip{\!\sm\!\!}{#1}}
\newcommand{\varp}[1]{\var\!\left(#1\right)}

\newcommand{\vara}{\varp{\fa}}
\newcommand{\varb}{\varp{\fb}}
\newcommand{\varc}{\varp{\fc}}

\newcommand{\varSName}{\Gamma}
\newcommand{\varS}[1]{\varSName_{#1}}

\newcommand{\varSB}[1]{\varSName^0_{#1}}
\newcommand{\varSBj}{\varSB{j}}

\newcommand{\fFun}{\bfm{f}}
\newcommand{\fFunp}[1]{\fFun\left(#1\right)}
\newcommand{\fFunj}[1]{\fFun_{#1}}

\newcommand{\tfFun}{\widetilde{\bfm{f}}}
\newcommand{\tfFunp}[1]{\tfFun\left(#1\right)}
\newcommand{\tfFunj}[1]{\tfFun_{#1}}

\newcommand{\tfun}{\widetilde{f}}

\newcommand{\tfunpp}[2]{\tfun_{#1,#2}}
\newcommand{\tfa}{\tfunpp{\bfm{a}}{j}}
\newcommand{\tfb}{\tfunpp{\bfm{b}}{j}}
\newcommand{\tfc}{\tfunpp{\bfm{c}}{j}}
\newcommand{\tfap}[1]{\tfa\left( #1 \right)}
\newcommand{\tfbp}[1]{\tfb\left( #1 \right)}

\newcommand{\dPijS}[2]{\upsilon_{#1,#2}}

\newcommand{\bigoo}[2]{\mathcal{O}_2\left(#1,#2\right)}

\newcommand{\dt}[4]{\begin{Vmatrix} #1 & #2 \\ #3 & #4 \end{Vmatrix}}
\newcommand{\dtP}[2]{\begin{Vmatrix} #1&#2\end{Vmatrix}}
\newcommand{\vtwo}[2]{\begin{pmatrix} #1 \cr #2 \end{pmatrix}}
\newcommand{\pc}{\bfm{p}}

\newcommand{\pcn}{\pc_n}
\newcommand{\pcnp}[1]{\pcn\left(#1\right)}

\newcommand{\T}{\bfm{P}}
\newcommand{\dT}{\Delta\bfm{P}}
\newcommand{\Tp}[1]{{\T}_{#1}}
\newcommand{\dTp}[1]{\dT_{#1}}
\newcommand{\vdm}{V}
\newcommand{\vdmp}[1]{\vdm\!\!\left( #1 \right)}
\newcommand{\vdmc}[2]{\vdm_{\bfm{t},\bfm{c},\left(#1,\dots,#2\right)}}

\newcommand{\qp}{\bfm{Q}}
\newcommand{\qpp}[1]{\qp_{#1}}

\newcommand{\starbase}{\sigma}
\newcommand{\lstar}{\starbase^*}
\newcommand{\lstarp}[1]{\lstar\left(#1\right)}
\newcommand{\ldstar}{\starbase^{**}}

\newcommand{\nrd}[1]{\| #1 \|}

\newcommand{\const}{\textrm{const}}

\newcommand{\closureB}[1]{\textrm{cl}\big(#1\big)}
\def\bfm#1{\boldsymbol{#1}}

\newcommand{\qqA}[1]{\mathcal{#1}^{\raise0.02ex\hbox{\scriptsize *}}}

\newcommand{\dist}[2]{\textrm{dist}\left( #1, #2 \right)}

\newcommand{\ie}{\textit{i}.{e}., }

\newtheorem{thm}{Theorem}
\newtheorem{cor}{Corollary}
\newtheorem{lem}{Lemma}
\newdefinition{rmk}{Remark}
\newdefinition{mydef}{Definition}
\newcommand{\lintrans}{A}
\newcommand{\lintransp}[1]{\lintrans_#1}

\newcommand{\asstrans}{B}
\newcommand{\asstransp}[1]{\asstrans_#1}

\newcommand{\asstranspj}{\asstransp{j}}
\newcommand{\assvec}{\bfm{v}}
\newcommand{\assvecp}[1]{\assvec_#1}
\newcommand{\assvecpj}{\assvecp{j}}
\newproof{pf}{Proof}
\newcommand{\mm}{m}
\newcommand{\mms}{\Sigma}
\newcommand{\mmsfp}[1]{\mms_{\backslash\left\{#1\right\}}}
\newcommand{\mmR}{\RR^{\mm}}

\newcommand{\mmUp}[1]{\bfm{U}_#1}
\newcommand{\mmF}{\bfm{g}}
\newcommand{\mmFt}{\widetilde{\bfm{\mmF}}}
\newcommand{\mmFp}[1]{\mmF\left(#1\right)}

\newcommand{\mmf}{g}
\newcommand{\mmfp}[1]{\mmf_#1}
\newcommand{\mmfpp}[2]{\mmfp{#1}\left(#2\right)}
\newcommand{\mmft}{\widetilde{g}}
\newcommand{\mmfpt}[1]{\mmft_#1}
\newcommand{\mmfppt}[2]{\mmfpt{#1}\left(#2\right)}

\newcommand{\mmM}{\mathcal{M}}
\newcommand{\coepsZero}{\ensuremath{\varepsilon}}

\newcommand{\coepsZerop}[1]{\coepsZero\left(#1\right)}
\newcommand{\coepsZeroS}{\ensuremath{\coepsZero^*}}

\newcommand{\coepsZeroiS}[1]{\ensuremath{\coepsZeroS_{#1}}}

\newcommand{\MNamep}[1]{\mathcal{M}_{\left(#1\right)}}
\newcommand{\Mi}{\MNamep{i}}
\newcommand{\Mii}{\MNamep{ii}}
\newcommand{\Miii}{\MNamep{iii}}
\begin{document}
\begin{frontmatter}



\title{On Planar Polynomial Geometric Interpolation}
 \author{Jernej Kozak\fnref{label2}}
 \ead{jernej.kozak@fmf.uni-lj.si}
 \ead[url]{https://arhiv-www.fmf.uni-lj.si/~kozak/wikiang/\\index.php?title=Main_Page}
\fntext[label2]{The author acknowledges the financial support from the Slovenian Research Agency (research core funding No.P1-0291)}



\begin{abstract}
In the paper, the planar polynomial geometric interpolation of data points is revisited. Simple sufficient geometric conditions that imply the existence of the interpolant are derived in general. They require data points to be convex in a certain discrete sense. Since the geometric interpolation is based precisely on the known data only, one may consider it as the parametric counterpart to the polynomial function interpolation. The established result confirms the H\"{o}llig-Koch conjecture on the existence and the approximation order in the planar case for parametric polynomial curves of any degree stated quite a while ago.
\end{abstract}




\end{frontmatter}


 
\section{Introduction}   \label{sec:introduction}

\begin{figure}[htb]
\hskip 0.3cm
  \begin{minipage}{0.45\textwidth}
		\includegraphics[width=\textwidth]{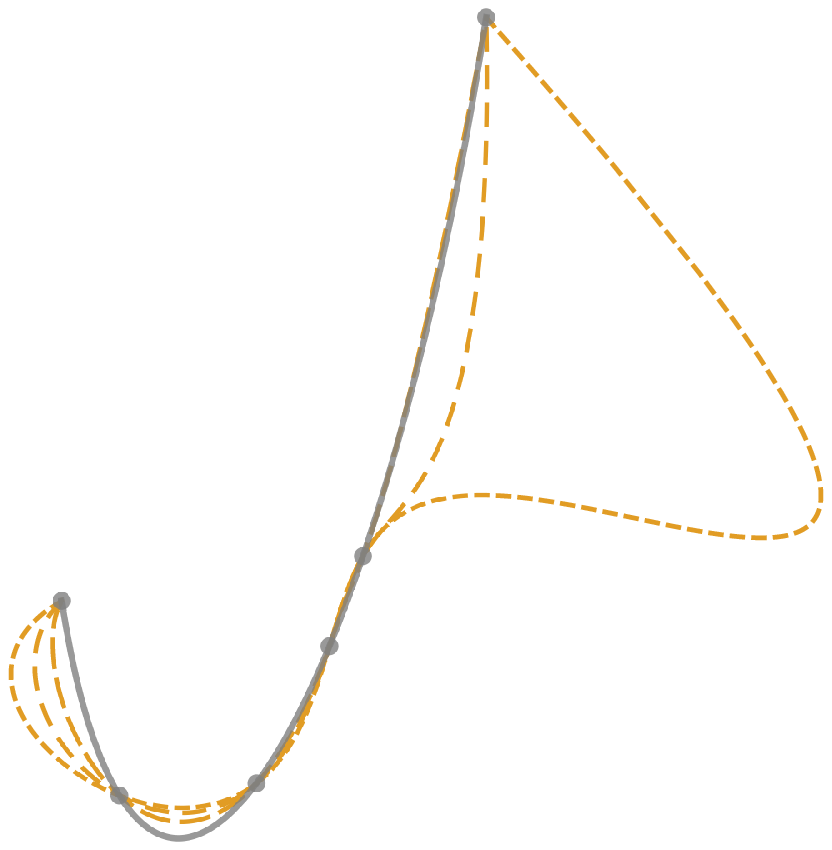}
  \end{minipage}
\hskip 0.5cm
 \begin{minipage}{0.45\textwidth}
 \centering 
	\includegraphics[width=\textwidth]{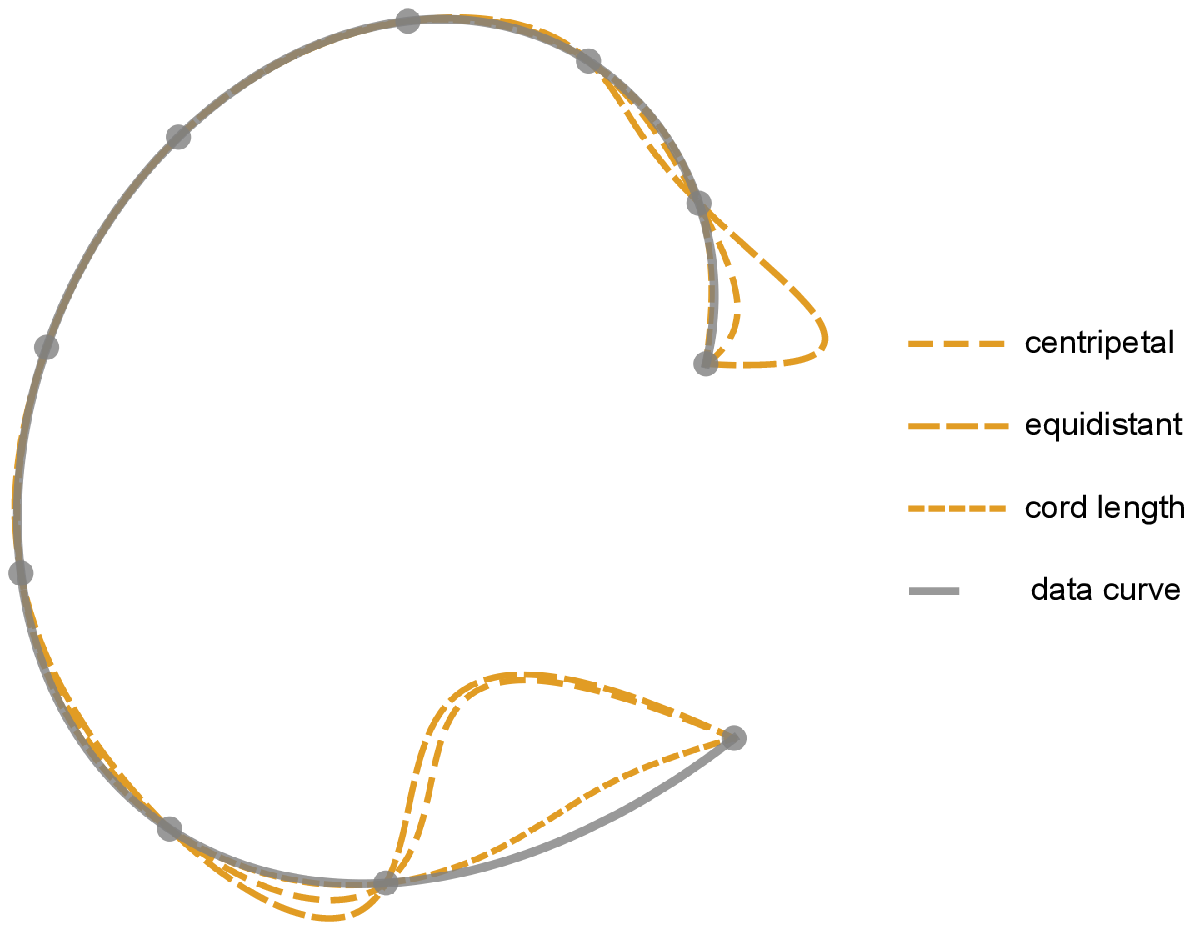}
  \end{minipage}
	\hskip 0.3cm
\caption{Examples of data   interpolated by a polynomial  parametric curve with parametrization chosen in advance.  Both data curves are smooth and convex.} 
\label{fig:figlin}
\end{figure}

The study of geometric polynomial interpolation problems was initiated in \citep{deBoor-Hoellig-Sabin-87-High-Accuracy} where Hermite cubic interpolation of two points, tangent directions, and curvatures was analyzed. It soon became apparent that geometric interpolation brings two important practical advantages
compared to linear polynomial interpolation schemes adapted to the curve case.
First, the interpolating curve depends on geometric quantities only:   data points, tangent directions, curvatures, \etc. 
As to the second one, we avoid the choice of parametrization in advance. This allows us to stick to lower degree polynomials, but use the implicit parametrization to keep the approximation power of the higher degree component-wise linear interpolants. In practical applications, the shape of an interpolating curve is often an issue. One would require that the interpolant at least roughly preserves the shape of the data it is based upon. However, a priori choice of parametrization combined with higher degree polynomials may make this goal hard to achieve. As an example, consider the points, sampled from smooth convex curves \figref{fig:figlin}. With three commonly used parametrizations applied, none of the interpolants is satisfying. But the same data interpolated geometrically, determine quite differently a satisfactory curve approximation \figref{figGeomInt}. So understandably a lot of work that followed \citep{deBoor-Hoellig-Sabin-87-High-Accuracy} investigated the geometric interpolation quite thoroughly. The problems considered included Lagrange interpolation of data points only, and its limit case, Hermite, 
where also data directions, data curvatures, \etc. are prescribed. 
If all data appear at one point, the Hermite case reduces to a Taylor one. 
But the problems involved are  generally nonlinear, so  the existence, the approximation order,  and the construction  are the burning issues.
The asymptotic approach  was often the tool to tackle the first two goals. There, the data is assumed to be taken from a sufficiently small part of a convex curve that could be locally parametrized by one of the components. 

In \citep{HolligKoch-1996}, a  natural conjecture on the geometric polynomial interpolation was made, based upon the authors'  previous work and results referenced therein. This conjecture in the special,  Lagrange planar case, predicts that  a polynomial curve of ${\textrm{degree}\ \le n}$ can if some appropriate assumptions are met interpolate $2n$ points, with optimal approximation order $2n$. The conjecture was confirmed in some particular cases, mostly for low degree polynomials, and in the asymptotic sense only 
 \citep[][\etal]{deBoor-Hoellig-Sabin-87-High-Accuracy, HolligKoch-1995, HolligKoch-1996,  MorkenScherer-1997, Scherer-2000, KozakZagar-2004, JaklicKozakKrajncZagar-2007b}. 

It is somewhat surprising how long the result for the general degree curves was slipping one's hands. In this paper, we close this gap for the planar Lagrange interpolation problem, and we present quite general geometric conditions that assure both the existence and the optimal approximation order for arbitrary degree cases. Recently,  a similar, positive result for the Taylor case polynomial curve expansion was given in \citep{Brysiewicz-2021} for quite a large infinite subset of $\NN$.  

The result of the paper is based  on an old  but far-reaching achievement  developed to  the  Bolzano–Poincaré–Miranda theorem for simplices in \citep[][ Theorem 2.1]{Vrahatis-2016}: if a set of  continuous  functions satisfies certain sign conditions at the simplex boundary, it must vanish at some interior point of the simplex.

The outline of the paper is  as follows. In Section~\ref{sec:Lagrange interpolation problem}, the Lagrange interpolation problem is defined, and the main theorems are listed. In Section~\ref{sec:Equations of the problem}, the functions and the corresponding equations  to be solved are derived. Section~\ref{sec:Analysis of polynomials}, the most technical part, contains the analysis of the functions involved. Based upon Section~\ref{sec:Analysis of polynomials} the last two sections provide the proofs of the theorems given in Section~\ref{sec:Lagrange interpolation problem}.

\section{Lagrange interpolation problem}\label{sec:Lagrange interpolation problem}

\begin{figure}[htb]
\hskip 0.3cm
  \begin{minipage}{0.40\textwidth}
		\includegraphics[width=\textwidth]{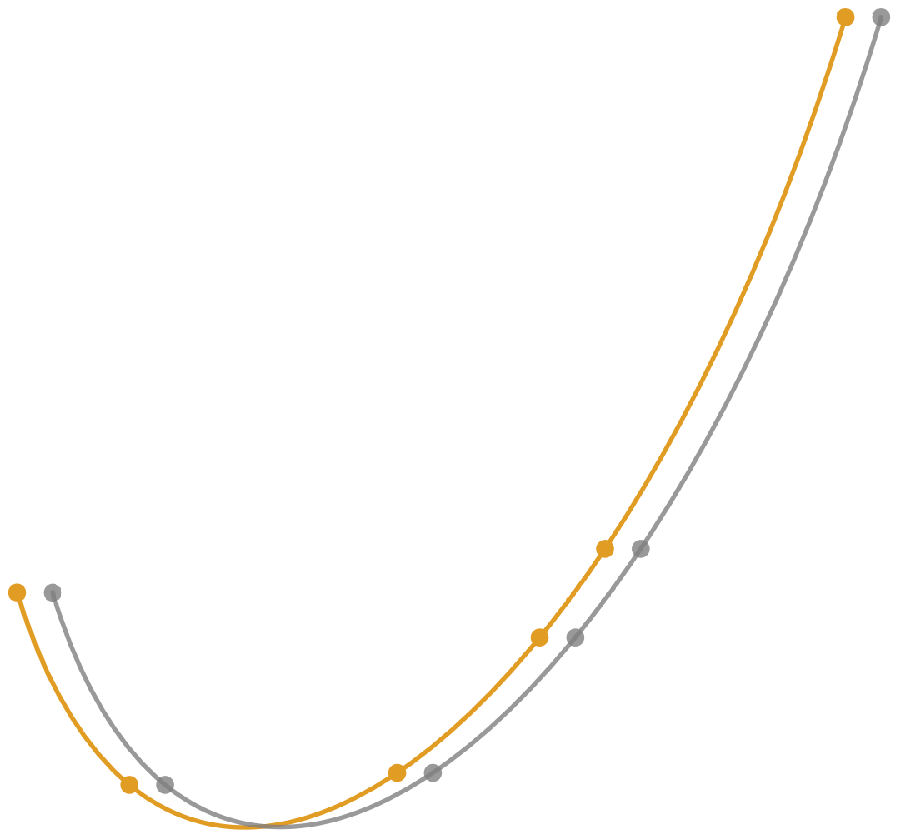}
  \end{minipage}
\hskip 1cm
  \begin{minipage}{0.40\textwidth}
	\includegraphics[width=\textwidth]{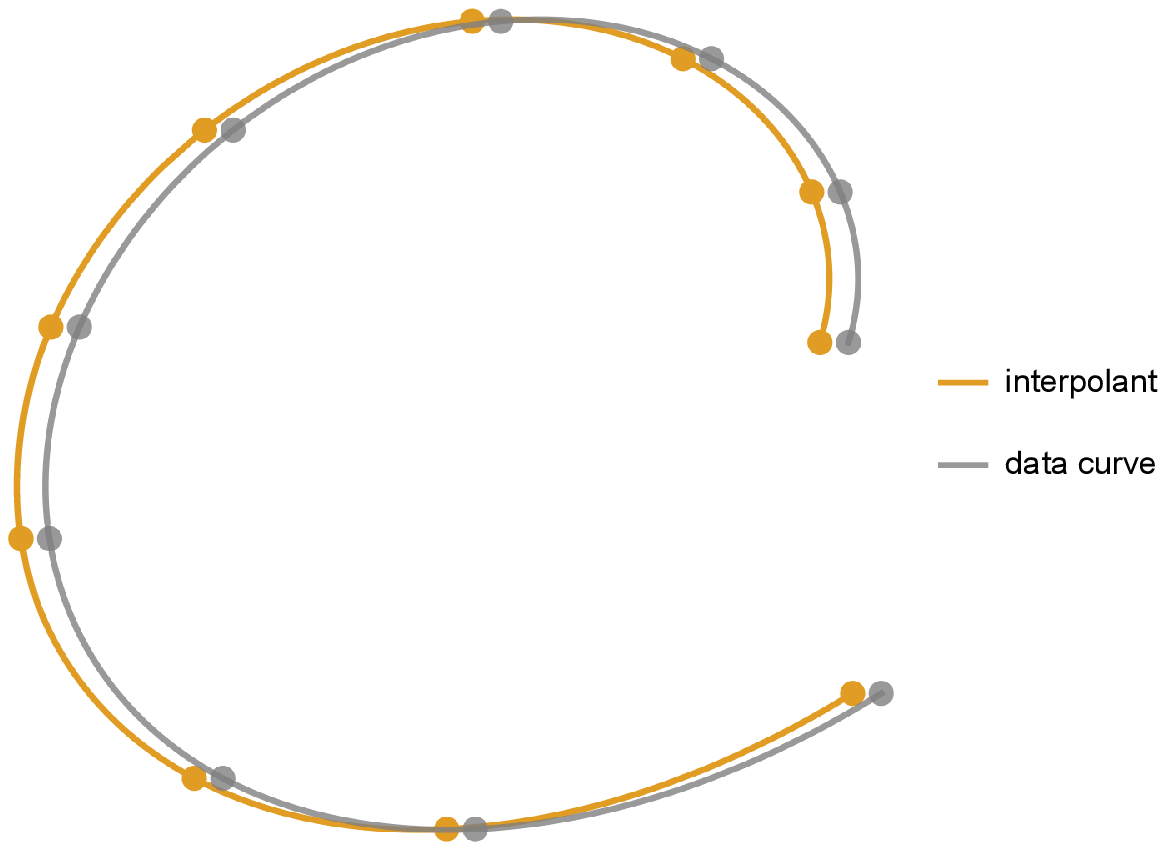}
  \end{minipage}
		\hskip 2.0cm
\caption{Data that satisfy \thmnref{theorem1}  interpolated by a cubic parametric curve (left), and data that satisfy \thmnref{theorem2} and corresponding polynomial interpolant of  degree  $5$ (right). Data curves are plotted with a small offset to the right since otherwise, any distinction between the original curve and its interpolant would not be visible.} 
\label{figGeomInt}
\end{figure}

The Lagrange interpolation problem is formulated as follows. 
Suppose that a sequence of $2n$ distinct points 
\begin{equation}\label{datadef}
\Tp{\ell} = \vtwo{a_\ell}{b_\ell} \in \RR^2, \   \ell = 0,1,\dots,2n-1,\quad \vv{a} := \vcc{a_\ell}{\ell=0}{2n-1},\  \vv{b} := \vcc{b_\ell}{\ell=0}{2n-1},
\end{equation}
is given. Find a 
parametric polynomial curve 
$$
  \pcn: [0,1] \to \RR^2
$$ 
of degree $\le n$ that
interpolates the given data points \eqref{datadef} at some parameter
values ${\ti{\ell} \in [0,1]}$ in increasing order, \ie
\begin{equation}\label{interpol-conditions}
\pcnp{t_\ell} = \Tp{\ell}, \quad \ell=0,1,\dots,2n-1.
\end{equation}
Since a linear reparametrization 
preserves the degree of a parametric polynomial curve, one can assume 
$\ti{0}:=0$ and $\ti{2n-1}:=1$, but the remaining parameters
\begin{equation*}
   \vt := \left(\ti{\ell} \right)_{\ell=1}^{2n-2}
\end{equation*}
are unknown, ordered as
\begin{equation*}
  \ti{0} = 0 < \ti{1} < \dots < \ti{2n-2} < \ti{2n-1} = 1.
\end{equation*}
Let $\Delta$ denote the forward difference, 
\ie
$$
\Delta w_\ell := w_{\ell+1} - w_{\ell}, \quad \dTp{\ell} := \Tp{\ell+1} - \Tp{\ell} .
$$
We have the following assertion.
\begin{thm}\label{theorem1} If, 
possibly after a linear 
transformation,
the
differences 
\begin{gather}\label{eqthmsupx} 
\deap{\ell} , \ \debp{\ell}, \quad  \ell = 0,1, \dots,2n-2, 
\end{gather}
 of the data \eqref{datadef} are all positive or all negative  and   the determinants
\begin{gather}\label{eqthmsupz}
\dtP{\dTp{\ell-1}}{\!\!\!\dTp{\ell}} =
\deap{\ell-1} \debp{\ell } -  \deap{\ell} \debp{\ell-1},
\quad \ell =1,2, \dots,2n-2, 
\end{gather}
are all of the same sign, the interpolation problem \eqref{interpol-conditions} has at least one  solution.
\end{thm}
Some obvious remarks are right on the spot. Even if the data \eqref{datadef} don't satisfy the assumptions  \eqref{eqthmsupx} and \eqref{eqthmsupz}, they may be satisfied after a proper linear data transformation. In particular, the role of both data components may be reversed, any data rotation or reflection is allowed, \etc.  All of the asymptotic approaches to the interpolation problem are covered here too. Generally, \thmnref{theorem1} gives the positive answer for all  data  that could be transformed to a sampled monotone convex  function curve by an affine transformation, as is the data case on \fignref{figGeomInt} (left). But the right side of the same figure requires more general assumptions.
\begin{thm}\label{theorem2} If the determinants
\begin{equation}\label{eqthmsupa}
\dtP{\dTp{j-1+k}}{\!\!\dTp{j-1+\ell}}
,\quad j =  1, 2, \dots,n-1; \, 0\le k < \ell \le n,
\end{equation}
are all of the same sign, the interpolation problem \eqref{interpol-conditions} has at least one solution.
\end{thm}
Note that the assumptions \eqref{eqthmsupa} are affinely invariant. 
The proof of \thmnref{theorem1} will carry most of the burden, and \thmnref{theorem2} will be confirmed with its help.
There is little hope that  simple necessary and sufficient conditions that guarantee existence could be found in general. In the case $n = 2$ the geometric interpolation problem has a solution iff \eqref{eqthmsupa} are satisfied as observed already in \citep{LachanceSchwartz-1991} and further elaborated in \citep{Morken-1995}. This case was generalized, together with necessary and sufficient conditions for the existence of the interpolant,  to any dimension in \citep{KozakZagar-2004} where  a particular form of \eqref{eqthmsupa} already appeared. But if the degree $n$ is larger, \ie $n \ge 3$ in the planar case, the analysis of all possibilities  becomes very complicated even for  planar cubic interpolation 
\citep{KozakKrajnc-2007, Krajnc-2009}.
If the data are not convex in the discrete sense of either one of the theorems, the existence of the solution depends heavily on the particular data set.
One would  expect that the results carry over to the Hermite interpolation problems, but we leave the confirmation to the future work.

\section{Equations of the problem}\label{sec:Equations of the problem}
\begin{figure}[htb]
	\centering
  \begin{minipage}{0.45\textwidth}
    \centering \includegraphics[width=\textwidth]{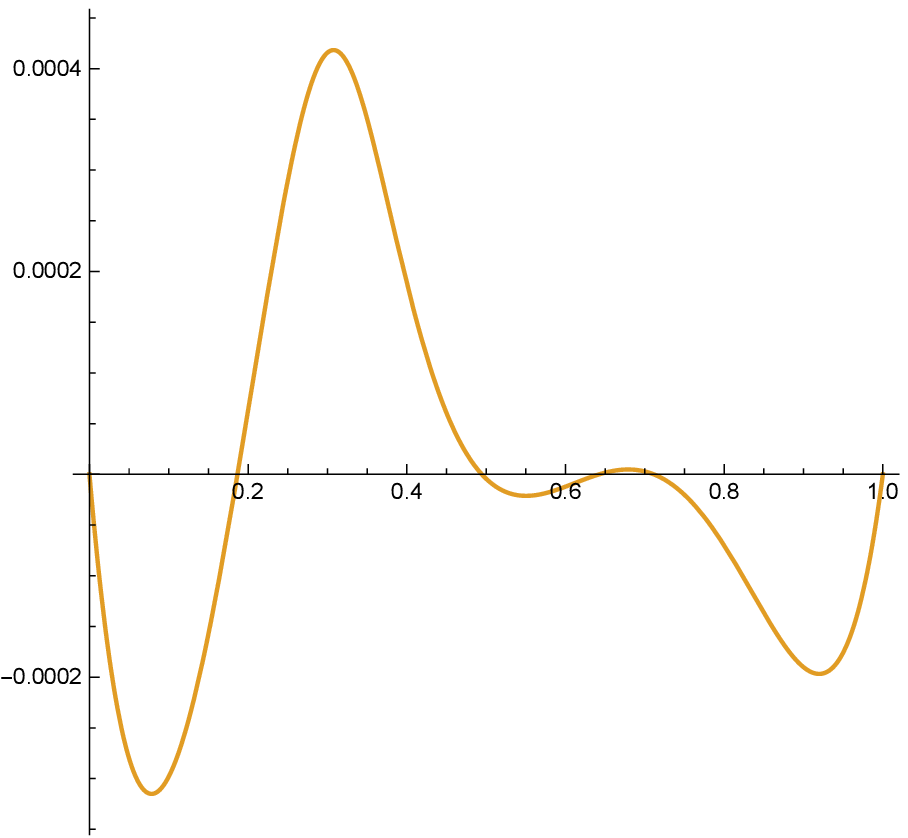}
  \end{minipage}
	\hfill
  \begin{minipage}{0.45\textwidth}
   \centering \includegraphics[width=\textwidth]{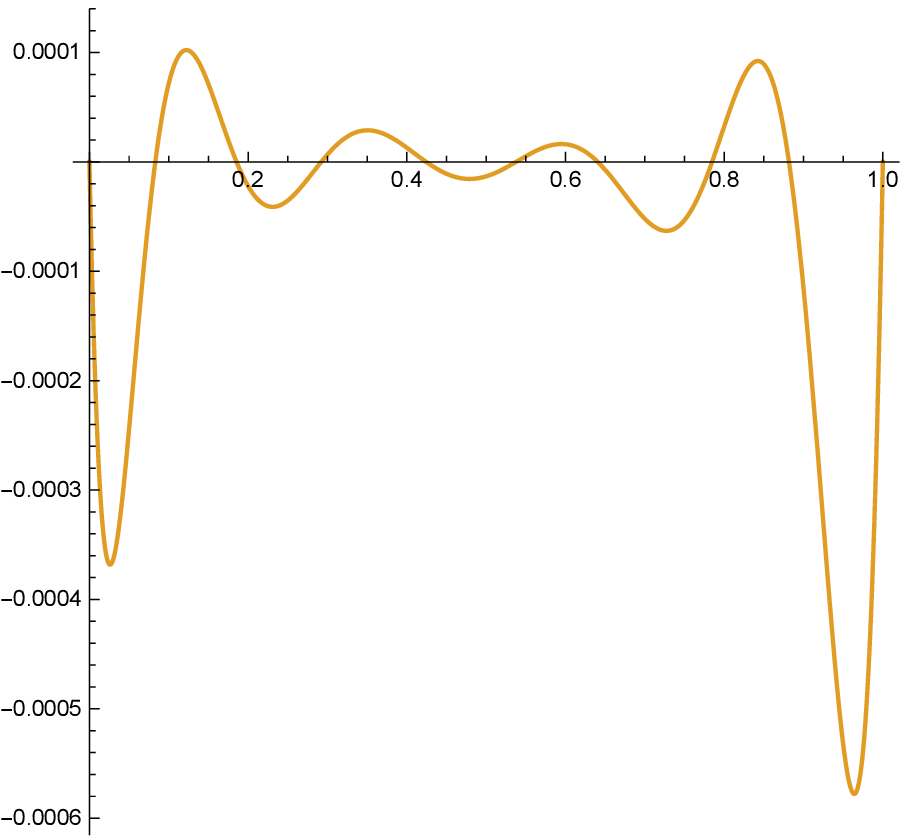}
  \end{minipage}
\caption{Normal reparametrization $\varphi$ and corresponding parametric difference graphs $\vv{g} \circ \varphi - \vv{p}$ of the examples outlined in  \fignref{figGeomInt}. } 
\label{figNorDiff}
\end{figure}

\begin{figure}[htb]
	\centering
  \begin{minipage}{0.45\textwidth}
    \centering \includegraphics[width=\textwidth]{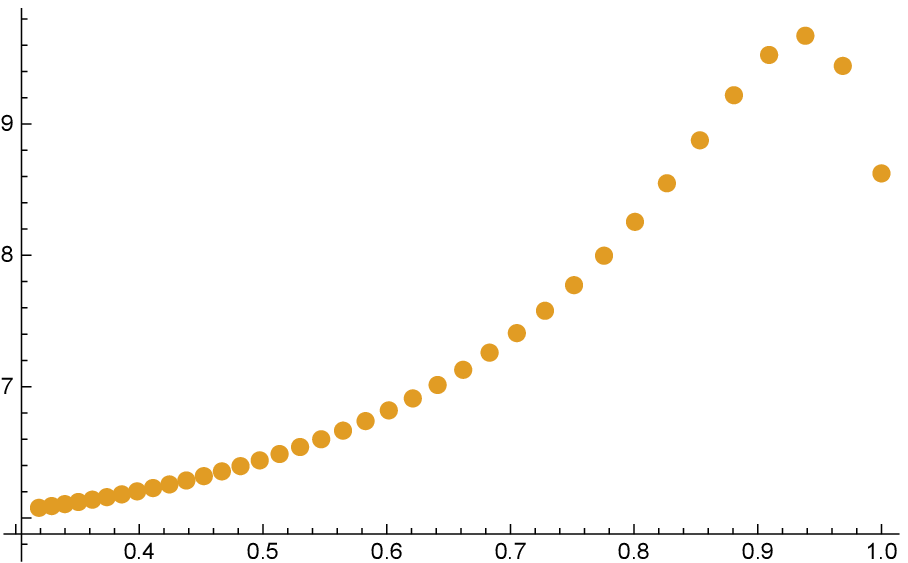}
  \end{minipage}
	\hfill
  \begin{minipage}{0.45\textwidth}
   \centering \includegraphics[width=\textwidth]{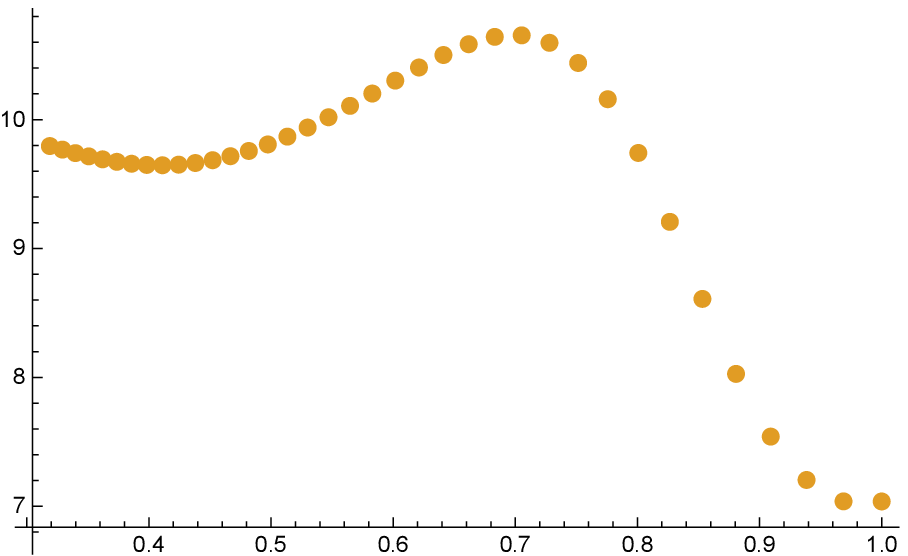}
  \end{minipage}
\caption{Numerical estimate of the approximation order for the data of  \fignref{figGeomInt}, left and right respectively. The abscissa denotes the  length of the data curve parameter relative to the length of the initial one.} 
\label{figApprOrd}
\end{figure}

As a vector in $\RR^{2n-2}$, the unknowns $\vt$  should belong to the   interior $\interior{\sm}$ of a ${(2n-2)}$-simplex $\sm$,
\begin{equation*}
\sm := \left\{ \vt \in \RR^{2n-2} | 0 \le \ti{1} \le \dots \le \ti{2n-2} \le 1\right\}.
\end{equation*}
The system of equations \eqref{interpol-conditions} has to determine 
the unknown
$\pcn$ as well as the parameters $\vt$. But the two tasks
can be separated if one can provide enough
linearly independent functionals, depending  on $\vt$ only, 
that map $\pcn$ to zero.
Divided differences, based upon $\ge~n+2$ values, are a natural 
choice.
Let us apply the divided differences
\begin{equation}\label{DD}
[t_{j-1}, t_{j}, \dots,t_{n + j}], \quad  j = 1,2,\dots,n-1,
\end{equation}
to both sides of \eqref{interpol-conditions}. 
Since ${\deg \ \pcn \leq n}$, the left side vanishes, 
and so should the right one. But the parameters $\ti{\ell}$ are distinct and
this condition becomes
\begin{equation}\label{system-closed-form}
   \displaystyle{\sum_{\ell=j-1}^{n+j}\frac{1}
  {\prod\limits_{\genfrac{}{}{0pt}{}{m=j-1}{m\ne
   \ell}}^{n+j}(\ti{\ell}-t_m)}}  
   \  \T_\ell= 0, 
   \quad
   j = 1,2, \dots, n-1.
\end{equation}
This nonlinear system depends on the data $\T_\ell$ and the unknowns $\vt$ only. 
If the solution ${\vt \in \interior{\sm}}$ is found, one can 
apply any algorithm inherited from the polynomial function interpolation to  construct  the interpolatory curve $\pcn$  component-wise.
For each~$j$ the system \eqref{system-closed-form} provides two equations based upon the first and the second component of the data. For our purpose, it will be more convenient to rewrite \eqref{system-closed-form} in a polynomial form. However, from a computational point of view, the equations \eqref{system-closed-form} are less sensitive to numerical cancellation, and they produced the examples of  \fignref{figGeomInt} in a couple of Newton steps starting with an equidistant initial guess.
Let $\vdmp{t_{j_{1}}, \dots, t_{j_{r}} }$  be a Vandermonde determinant based upon $t_{j_{1}}, \dots, t_{j_{r}}$,
\begin{equation}\label{vandmdef}
\vdmp{t_{j_{1}}, \dots, t_{j_{r}} }  :=  \deter{ t_{j_{\ell}}^{i-1} }{i=1;\ell=1}{r;r} =  \prod_{ \ell = 1}^{ r} \prod_{i=1}^{\ell-1} \left( t_{j_\ell} - t_{j_i} \right),
\end{equation}
and let
\begin{equation*}
\vdmp{t_{j_{1}}, \dots,t_{j_{r}};c_{j_{1}},  \dots,c_{j_{r}} } 
\end{equation*}
denote a determinant, obtained from \eqref{vandmdef} by replacing the last determinant row \newline
${\left( t_{j_{1}}^{r-1},  \dots,t_{j_{r}}^{r-1} \right)}$ by the values $\left( c_{j_{1}},  \dots,c_{j_{r}} \right)$,
\begin{equation*}
\vdmp{t_{j_{1}},  \dots,t_{j_{r}};c_{j_{1}},  \dots, c_{j_{r}} }  
= \sum_{\ell=1}^r \left( -1\right)^{r+\ell} c_{j_\ell} \vdmp{t_{j_{1}}, \dots, t_{j_{\ell-1}}, t_{j_{\ell+1}},\dots,t_{j_{r}} }.
\end{equation*}
The divided difference \eqref{DD} applied to a sequence $\vcc{c_\ell}{\ell = j-1}{\ell=n+j}$ as a quotient of the Vandermonde determinants reads
\begin{equation*}\label{eq:quocvdm}
\frac{
\vdmp{t_{j-1},  \dots,t_{n+j};c_{j-1},  \dots, c_{n+j} }
}{\vdmp{t_{j-1}, \dots,t_{n+j} } }.
\end{equation*}
%
Thus, if we multiply \eqref{system-closed-form} by $\vdmp{t_{j-1}, \dots,t_{n+j} }$ we obtain the polynomial form of the system
\begin{align}\label{system-pol-form}
& \fap{t_{j-1},  \dots,t_{n+j}}   = 0, \nonumber \\[-4pt]
 & \quad\quad\quad\quad\quad\quad\quad\quad\quad\quad\quad  j = 1,2,\dots,n-1,\\[-4pt]
    & \fbp{t_{j-1},  \dots,t_{n+j}}   = 0, \nonumber 
\end{align}
where 
\begin{equation*}
\fcp{t_{j-1},  \dots,t_{n+j}} 
:=  \vdmp{t_{j-1},  \dots,t_{n+j};c_{j-1},  \dots,c_{n+j} } .
\end{equation*}
Here and throughout the paper,  coefficients $\vv{c}$ will be used as a placeholder for either $\vv{a}$ or $\vv{b}$ meaning that any relation or text involving $\vv{c}$ holds for either of them. Note that $\fc$ depends linearly on the constants 
$
c_{j-1}, c_j, \dots, c_{n+j}
$
 only. It is also invariant under data translation since divided differences map constants to zero. So $\fc$ depends linearly on data differences $\decp{\ell}$,
\begin{align}\label{fdifdef}
& \fcp{t_{j-1},  \dots,t_{n+j}}  =  \nonumber \\
& \quad = \sum_{k=j}^{n+j} \left( -1\right)^{n+j-k} (c_{k} - c_{j-1})\vdmp{t_{j-{1}}, \dots, t_{{k-1}}, t_{{k+1}},\dots,t_{n+j} } \\
& \quad  = \sum_{r=j}^{n+j} \decp{r-1} \sum_{k=r}^{n+j} \left( -1\right)^{n+j-k}\vdmp{t_{j-{1}}, \dots, t_{{k-1}}, t_{{k+1}},\dots,t_{n+j} }. \nonumber
\end{align}
This shows that the equations derived follow the affine nature of the interpolation problem,
and it explains the following remark.
\begin{rmk}\label{rem:remark1}
 For each $j$ separately, $ 1\le j \le n-1$, we may replace  the data 
 \begin{equation}\label{remdat}
\begin{pmatrix} 
\dTp{j-1} & \dTp{j}& \dots &\dTp{n+j-1}
\end{pmatrix}
\end{equation}
which in \eqref{system-pol-form}  determine the function pairs $\fa$, and $\fb$,
by 
\begin{equation*}\label{remdata}
\begin{pmatrix} 
\lintransp{j}\,\dTp{j-1} & \lintransp{j}\,\dTp{j}& \dots &\lintransp{j}\,\dTp{n+j-1}
\end{pmatrix},
\end{equation*}
where $\lintransp{j}: \RR^2 \to \RR^2$ is any  nonsingular linear transformation. The system of equations generated by the modified functions is equivalent to the original one.
\end{rmk}
%
%
The following remark contributes an error term in case the data are sampled from a data curve.
\begin{rmk}\label{rem:remark2}
 Suppose that the data  
$$
 \Tp{\ell} = \vv{g}(\xi_{\ell}), \ \ell = 0,1,\dots,2n-1,
$$
are taken from a smooth parametric curve $\vv{g}: [\xi_0,\xi_{2n-1}] \to \RR^2$ at increasing parameter values
$
\xi_0 < \xi_1 < \dots < \xi_{2n-1}.
$
If the solution $\vv{t}\in \interior\sm$ of  \eqref{system-pol-form} exists,  the function case polynomial interpolation provides  the remainder  in the parametric one too. Let $\varphi: [0,1] \to  [\xi_0,\xi_{2n-1}]$ be any regular reparametrization of $\vv{g}$ that satisfies 
$$
\varphi(t_\ell) = \xi_\ell, \quad \ell=0,1,\dots,2n-1.
$$ 
Then
$$
\vv{g}(\varphi(t)) = \pcnp{t} + \left(\prod_{\ell=0}^{2n-1}\left(t-t_\ell\right)\right) \,  [\xi_{0}, \xi_{1}, \dots,\xi_{2n-1},\varphi(t)]\vv{g}, \quad t \in [0,1].
$$
As an example, the normal  reparametrization $\varphi$, introduced in  \citep{Degen-approximation-95}, is often a sensible choice \figref{figNorDiff}.
\end{rmk}
\begin{rmk}
If the data curve $\vv{g}$ in \remnref{rem:remark2} 
is additionally convex on $[\xi_a,\xi_b]$,
the points $ \Tp{\ell} $ satisfy the assumptions of  Theorem~\ref{theorem2}  for any 
increasing choice of parameters ${\xi_\ell \in [\xi_a,\xi_b]}$.
So by 
\citep[][Theorem~4.6]{JaklicKozakKrajncZagar-2007b} 
the interpolating polynomial
parametric curve $\vv{p_n}$  approximates $\vv{g}$ with the optimal approximation
order $2n$. \fignref{figApprOrd} shows a numerical estimate of the approximation order
for the data of \fignref{figGeomInt}  as a function of shrinking  data curve parameter interval computed as in 
\citep[][at several spots]{deBoor-Guide-2001}. 
The numerical evidence clearly confirms the expected result, ${2n = 6}$~(left), and ${2n=10}$~(right).
\end{rmk}

\section{Analysis of polynomials $\fc$}\label{sec:Analysis of polynomials}

\begin{figure}[htb]
	\centering
 \begin{minipage}{0.45\textwidth}
    \centering \includegraphics[width=\textwidth]{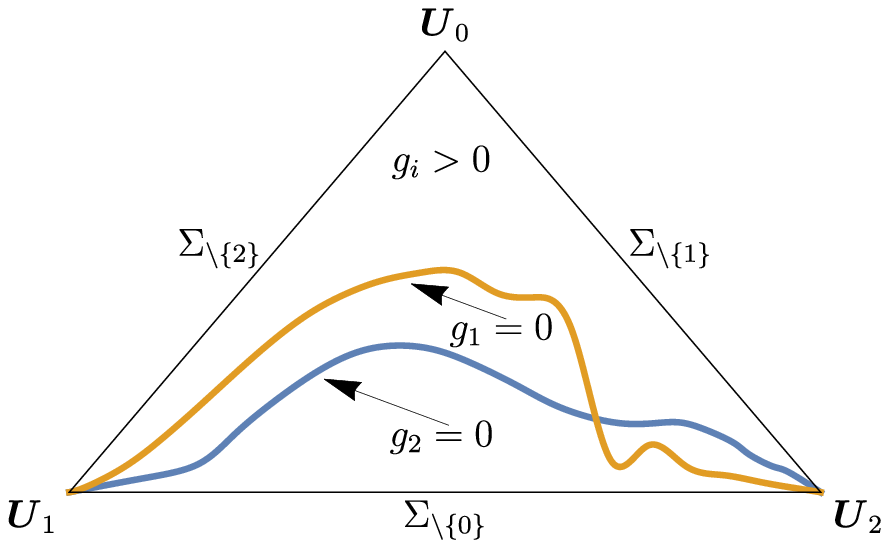}
 \end{minipage}
 \begin{minipage}{0.45\textwidth}
   \centering \includegraphics[width=\textwidth]{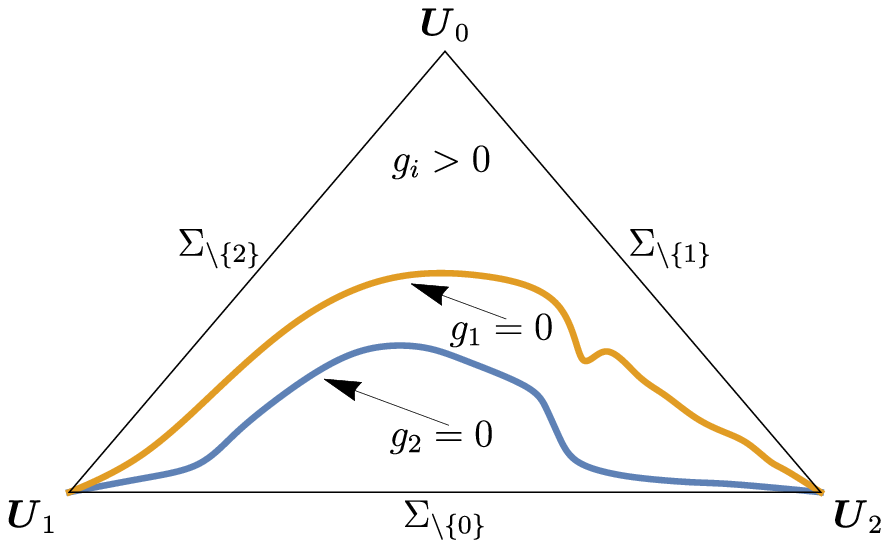}
  \end{minipage}
	\caption{The sign patterns of $\mmF$ at the boundary $\partial \mms$ of the left and the right image in the figure match, and the assumptions of Theorem~\ref{thmgBPM} are not satisfied in either case. }
\label{figgBPMa}
\end{figure}

Somewhat sloppy,  we shall consider the polynomial $\fc$ wherever needed also as a function with the domain extended naturally,
$$
\fc: \RR^{2n-2} \to \RR:\, \left( t_\ell \right)_{\ell=1}^{2n-2} \mapsto \fcp{t_{j-1},  \dots,t_{n+j}},
$$
and
\begin{equation*}
\fFunj{j} := \left( \funpp{\bfm{a}}{j}, \funpp{\bfm{b}}{j} \right)^T\!, \quad \fFun := \left( \funpp{\bfm{a}}{1}, \funpp{\bfm{b}}{1}, \dots,  \funpp{\bfm{a}}{n-1}, \funpp{\bfm{b}}{n-1} \right)^T. 
\end{equation*}
This helps us to define varieties $\varc$, $\varp{\fFunj{j}}$ and $\varp{\fFun}$, where the variety definition of a scalar function $g$,
$$
\varp{g} := \left\{ \vv{t} \in \RR^{2n-2}  \, \big| \, g\left(\vv{t}\right) = 0\right\},
$$
is naturally extended to the vector function case  $\varp{\bfm{g}}$.
The system \eqref{system-pol-form} is equivalent to \eqref{system-closed-form} except for possible extraneous solutions at the simplex boundary  
${\boudp{\sm} := \sm \setminus \smint}$, and there is to prove $ \varp{\fFun} \cap \interior{\sm} \neq \emptyset$.
Of particular interest will be a variety part 
\begin{equation*}
\varisp{g} := \closureB{\varp{g} \cap \interior{\sm} }
\end{equation*}
the restriction of $\varp{g}$ to $\sm$ that drops out the variety boundary points that can't be reached  from $\smint$.
The simplex $\sm$ requires some additional notation.  It is a convex hull of points,
$$
\sm = \cop{\qpp{0},\dots,\qpp{2n-2}},
$$
where
\begin{equation}\label{qidef}
\qpp{\ell} := ( \underbrace{0,\dots,0}_{\ell},  \underbrace{1,\dots,1}_{2n-2-\ell} )^T, \quad \ell= 0, 1,\dots,2n-2.
\end{equation}
The simplex faces will be denoted by
\begin{align*}
& \smfp{\ell_1,\dots,\ell_r} =  \cop{\qpp{j}}_{j \in \{0,\dots,2n-2 \} \backslash \{ \ell_1,\dots,\ell_r \}}.
\end{align*}
In this notation, the boundary of the simplex $\sm$ is determined by the {$(2n-3)\textrm{-simplex}$}  faces
\begin{equation}\label{eq:facef}
\smfp{\ell},  \quad \ell= 0, 1,\dots,2n-2.
\end{equation}
\begin{figure}[htb]
	\centering
  \begin{minipage}{0.48\textwidth}
    \centering \includegraphics[width=\textwidth]{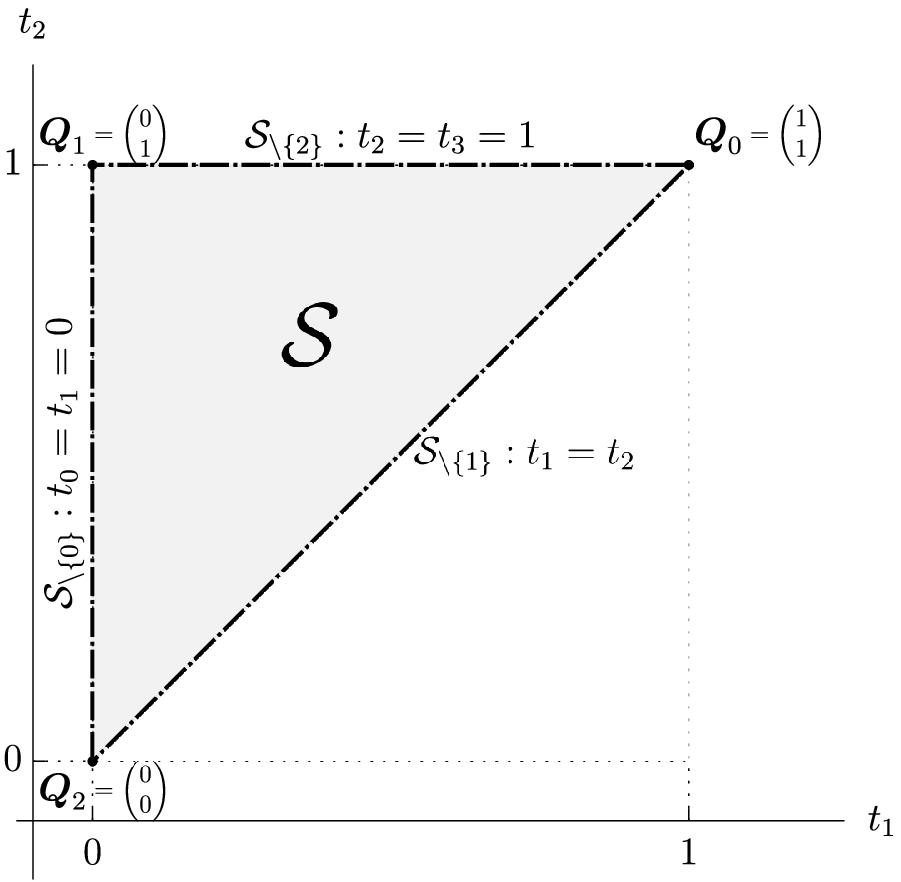}
  \end{minipage}
\caption{The 
case $n=2$: the simplex $\sm \subset \RR^2$, the points $\qpp{0},\qpp{1},\qpp{2}$, and the faces $\smfp{0},\smfp{1},\smfp{2}$. }
\label{figSdefL}
\end{figure}
The simplest case $n=2$ is shown in \fignref{figSdefL}. The first observation is straightforward.
Since precisely the components $\ell$ and ${\ell +1}$ of  
$\qpp{\ell}$ 
are by  \eqref{qidef} distinct, 
the $(2n-3)$-simplex $\smfp{\ell} $ is the maximal face of $\sm$ characterized by $t_\ell = t_{\ell+1}$.  Given $j$, $1\le j \le n - 1$, suppose that 
$$
j-1 \le r_1, r_2  \le n + j -1, \  r_2 \neq r_1, \  \  t_{r_2} = t_{r_1}.
$$
We then deduce  from \eqref{fdifdef} 
\begin{align}\label{fcval}
 & \fcp{t_{j-1}, \dots,  t_{r_2-1},t_{r_1},t_{r_2+1},\dots,t_{n+j}} =  \nonumber \\[-4pt]
 & \\[-4pt]
& \quad\quad \left( -1\right)^{n+r_2 -j} \left(  c_{r_2}   - c_{r_1} \right) \vdmp{t_{j-1}, \dots, t_{r_2-1}, t_{r_2+1},\dots,t_{n+j} } . \nonumber
\end{align}
Let us define  the signum function  by
\begin{equation*}
\sign \ w :=
 \left\{
\begin{array}{ll}
\ \ \ 1, & w > 0, \\
\ \ \ 0, & w = 0, \\
-1,& w < 0,\\
\end{array}
\right. 
\end{equation*}
and let us extend it  to a vector argument component-wise, ${\sign\left(w_i \right) := \left( \sign\  w_i \right)}$.
The choice  $r_2 = r_1 + 1 = \ell + 1$ in \eqref{fcval} proves the following lemma.
\begin{lem}\label{lem:lemma1} 
Let $ j-1 \le \ell \le n+j-1$. Then
$$
\sign \ \fcp{\vv{t}} = \left(-1\right)^{n+j-1-\ell} \sign \ \Delta c_{\ell} , \quad \vv{t} \in \interior{\smfp{\ell} },
$$
and further, 
\begin{align}\label{eq:sbound}
 & \sign \ \fcp{\vv{t}}_{|_{\vv{t} \in  \smfp{ \ell,k}}} =\\%
& \ = \left\{
 \begin{array}{ll}
0, &  k = j-1,\dots,\ell-1,\ell+1,\dots,n+j-1, \\
 \sign \ \fcp{\vv{t}}_{|_{\vv{t} \in  \,\interior{\smfp{\ell}}}} , &  k = 0,\dots,j-2,n+j,\dots,2n-2.
\end{array}
\right. \nonumber
\end{align}
\end{lem}
So the sign of  $\fc$ at the open faces  
$$
\interior{\smfp{\ell} },\quad \ell = j-1,\dots,n+j-1,
$$
is known for each $j$, and this fact even simplifies if all differences $\Delta c_{\ell}$ are of the same sign. This paves the way to an application of the Bolzano–Poincaré–Miranda theorem on simplices 
\citep[][Theorem 2.1]{Vrahatis-2016}. 
%
%
For  the reader's convenience, we recall the theorem here.
\begin{thm} \label{thmgBPM}
(\citep[][Theorem 2.1]{Vrahatis-2016})
Assume that 
$
\mms := \,
\cop{\mmUp{0},\dots,\mmUp{m}}
$ 
is an $\mm$-simplex in $\mmR$ with vertices $\mmUp{\ell} \in \RR^\mm$. Let 
$
\mmF :=  \left( \mmf_1, \dots,\mmf_\mm \right)^T: \mms \to \RR^\mm
$
be a continuous function such that 
$
\mmfp{i}\left(\mmUp{\ell}\right) \ne 0, \, i = 1,2,\dots,\mm; \, \ell = 0,1,\dots,\mm,
$
and 
$
{\mmFp{\bfm{x}} \ne \bfm{0}, \ \bfm{x} \in \partial\mms.}
$
Assume that the vertices $\mmUp{\ell}$ are reordered such that the following hypotheses are fulfilled:
\begin{enumerate}[label=(\roman*)]\setlength\itemsep{-4pt}
\item $ \sign \, \mmfpp{i}{\mmUp{i}} \sign \, \mmfpp{i}{\bfm{x}} =  -1, \quad \bfm{x} \in \mmsfp{i}, \  i = 1,\dots,\mm $, 
\item $ \sign \, \mmFp{\mmUp{0}} \ne  \sign \, \mmFp{\bfm{x}},  \quad \bfm{x} \in \mmsfp{0}. $
\end{enumerate}
Then, there is at least one $\bfm{x} \in \interior{\mms}$ such that $\mmFp{\bfm{x}}  =    \bfm{0}$.
\end{thm}
The sign assumptions of \thmnref{thmgBPM} are too demanding for a straightforward application in our case since by \lemnref{lem:lemma1} functions involved  vanish at a significant part of the boundary $\partial\sm$. An intermediate step is needed. We illustrate it  with two simple examples shown in \fignref{figgBPMa}.
Both functions ${\mmF: \mms \to \RR^2}$ there have an equal sign pattern at the boundary $\partial \mms$,
and there are almost none of the required assumptions  fulfilled. So  \thmnref{thmgBPM} can't be applied. But there is an ${\bfm{x} \in \interior \mms}$  such that $\mmFp{\bfm{x}} = 0$ in \fignref{figgBPMa}~(left), and no such interior point exists in the right part  of the figure. 
This encourages one to seek a modification $\mmF \to \mmFt$, based upon  some additional information on $\mmF$ at or very close to the boundary, such that the theorem could be put to work on $\mmFt$ where appropriate. Of course, the modification must preserve the zero set
$
{\{ \bfm{x} \in \interior \mms| \, \mmFp{\bfm{x}} = \bfm{0} \}.}
$
The following observation offers a simple sufficient condition that helps us to verify this.
\begin{lem}\label{lemmodstep}
Suppose that $\mmF :=  \left( \mmf_1, \dots,\mmf_\mm \right)^T: \mms \to \RR^\mm$ is  a continuous function. Choose  $\mmfp{i}$ that is to be replaced by a continuous $\mmfpt{i}: \mms \to \RR$, and let
$$
\mmM := \{ \bfm{x} \in \interior \mms  | \, \left( \mmfpp{i}{\bfm{x}} = 0 \wedge \mmfppt{i}{\bfm{x}} \ne 0 \right) \vee \left( \mmfppt{i}{\bfm{x}} = 0 \wedge \mmfpp{i}{\bfm{x}} \ne 0 \right)\}.
$$
If there exists a component $\mmfp{\ell}$, \, $1 \le \ell \le m$, 
such that
\begin{align*}
& \{ \bfm{x} \in \interior \mms  | \, \mmfpp{\ell}{\bfm{x}} = 0 \} \cap \mmM = \emptyset,
\end{align*}
then the replacement $\mmfp{i} \to \mmfpt{i}$ doesn't alter the set
$
\{ \bfm{x} \in \interior \mms| \, \mmFp{\bfm{x}} = \bfm{0} \}.
$
\end{lem}
\begin{pf}
The set $\mmM \subset \interior \mms$ determines the points  where the zero sets of $\mmfp{i}$ and $\mmfpt{i}$ differ. Since $\mmfp{\ell}$ doesn't vanish there, so can't $\mmF$, and the replacement $\mmfp{i} \to \mmfpt{i}$ does not influence  the zero set of $\mmF$.
\qed
\end{pf}
\begin{figure}[htb]
	\centering
 \begin{minipage}{0.45\textwidth}
    \centering \includegraphics[width=\textwidth]{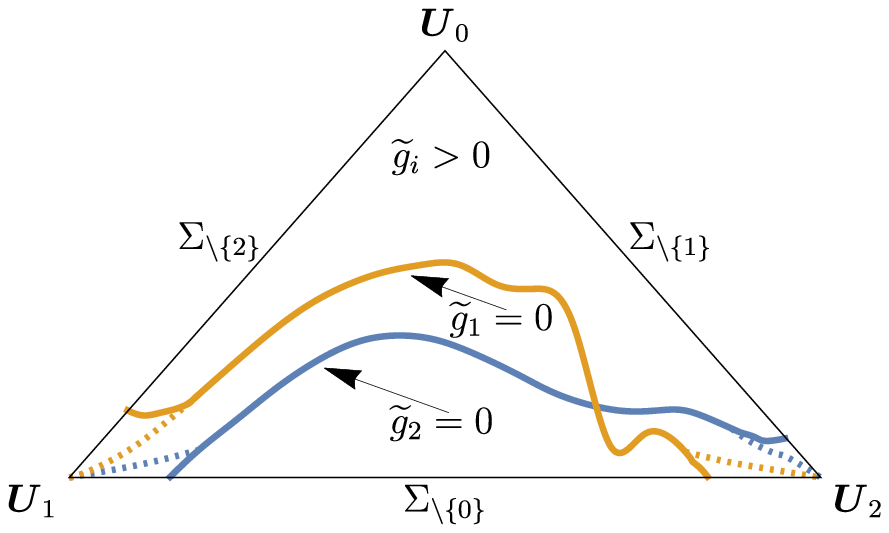}
 \end{minipage}
 \begin{minipage}{0.45\textwidth}
   \centering \includegraphics[width=\textwidth]{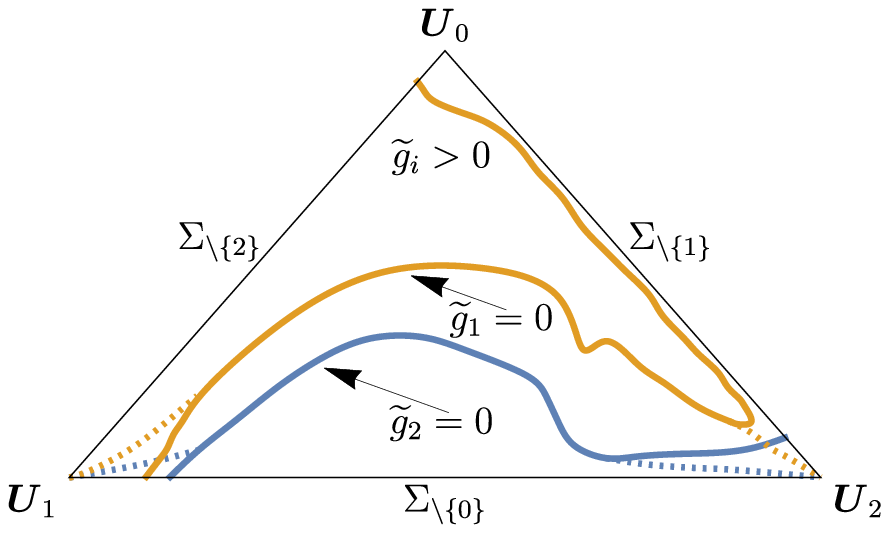}
  \end{minipage}
	\caption{ Modifications that satisfy \lemnref{lemmodstep} applied to the examples of  \fignref{figgBPMa}. The assumptions of \thmnref{thmgBPM} are met on the left, but not all on the right part of the figure.}
\label{figgBPMb}
\end{figure}
Recall \fignref{figgBPMa}. After several  steps  where each of them satisfies the assumptions of \lemnref{lemmodstep} the examples are modified as shown in \fignref{figgBPMb}. The function $\mmFt$ of the left image  fulfills the assumptions of \thmnref{thmgBPM} thoroughly, but for the right counterpart understandably this is not possible. All the other assumptions been satisfied, it still misses \emph{(ii)}.
%
%

To apply \thmnref{thmgBPM} a proper association between function $\fc$  and simplex faces $\smfp{\ell}$ has to be found. 
For each of the functions \eqref{system-pol-form} identified by the index pair $\left(\vv{c},j\right)$ we have to select  a  face $\smfp{\lstar}$, determined by a map
\begin{equation}\label{lstar}
 \left\{ \vv{a}, \vv{b} \right\} \times \left\{ 1,\dots, n-1 \right\} \to \left\{ 0,\dots,2 n - 2 \right\}: \, \left(\vv{c},j\right) \mapsto \lstar := \lstarp{{\vv{c},j}} 
\end{equation}
such that
\begin{equation}\label{pog1}
\sign \fcp{\vv{t}} \ \sign\fcp{
\qpp{\lstar}
} \ =  - 1, \quad \vv{t} \in \smfp{\lstar}.
\end{equation}
The map \eqref{lstar} should be injective, and combined with \eqref{pog1}  it represents \emph{(i)} of \thmnref{thmgBPM}.  
There are only $2n-2$ functions in \eqref{system-pol-form} and $2n-1$  simplex faces in \eqref{eq:facef}. One face is left, \ie $\smfp{\ldstar}$, at which the final  condition  \emph{(ii)} of \thmnref{thmgBPM} is imposed,
\begin{equation}\label{pog2}
\sign\fFunp{\bfm{t}} \ne \sign\fFunp{\qpp{\ldstar} },\quad \bfm{t} \in \smfp{\ldstar}.
\end{equation}
Quite clearly, given $j$ and $\fc$, the candidates for $\lstarp{{\vv{c},j}} $ are only the indices $j-1,j,\dots,n+j-1$. Note that \eqref{pog1} requires $\sign \fc$ to
be constant at the particular simplex, which partially collides  with \eqref{eq:sbound} established by Lemma~\ref{lem:lemma1}. So modifications of the functions $\fc$ are needed, 
and \lemnref{lemmodstep} together with the examples of \fignref{figgBPMa} will help us in the construction.

Let us elaborate \eqref{fcval} one step further. Note first that
\begin{align}\label{dvdm}
 & \frac{\partial}{\partial t_{j_i}} V\left( t_{j_1},\dots, t_{j_r}\right)_{\big|_{t_{j_i} = t_{j_m}}} = \nonumber \\[-4pt]
 & \\[-4pt]
& \quad\quad\quad = (-1)^{r-i} V\left( t_{j_1},\dots,t_{j_{i-1}},t_{j_{i+1}},\dots,t_{j_{r}}\right) \omj{1}{r}{i,m}{t_{j_m}} \nonumber .
\end{align}
From \eqref{fcval} and \eqref{dvdm} it is straightforward to derive first order terms of the Taylor expansion of $\fc $. Let 
$$
j - 1 \le r_1, r_2, r_3, r_4 \le n+j, \quad r_2 \ne r_1, r_4 \ne r_3, \ r_2 \le r_4,
$$
and
\begin{align}\label{lindef}
   &  \linp{r_1,r_2,r_3,r_4}  :=  \nonumber \\[-4pt]
	& \\[-4pt]
& \quad\quad 	\left(-1\right)^{r_4 - r_2}  
 \left( c_{r_4}   - c_{r_3},   c_{r_2}   - c_{r_1} \right) 
\begin{bmatrix}
	- \, \omdcp{r_2,r_4}{t_{r_1}} & 0 \\
	0 & \, \omdcp{r_2,r_4}{t_{r_3}}
\end{bmatrix}
\vdva{t_{r_2} - t_{r_1}}{t_{r_4} - t_{r_3}} , 
\nonumber
\end{align}
with
$$
\omdp{r_2,r_4}{z} := \om{j-1}{n+j}{r_2,r_4}{z}.
$$
Then $\fc$ at $t_{r_2} \approx t_{r_1}$ and $t_{r_4} \approx t_{r_3}$ expands as
\begin{align}\label{fcexp}
\fcp{t_{j-1}, \dots,t_{n+j}} & =  \vdmp{t_{j-1}, \dots, t_{r_2-1}, t_{r_2+1},\dots, \dots, t_{r_4-1}, t_{r_4+1},\dots,t_{n+j} }  \cdot \nonumber \\
&   \quad \cdot  \linp{r_1,r_2,r_3,r_4} + \bigoo{ t_{r_2} - t_{r_1}}{t_{r_4} - t_{r_3}}.
\end{align}
Here,
$
\bigoo{u}{v}
$
denotes the remainder involving terms  ${u^i v^{m-i}, 0 \le i \le m, m = 2, \dots, }$ with polynomial coefficients that depend on $\vt$.

\begin{figure}[htb]
	\centering
  \begin{minipage}{0.48\textwidth}
    \centering \includegraphics[width=\textwidth]{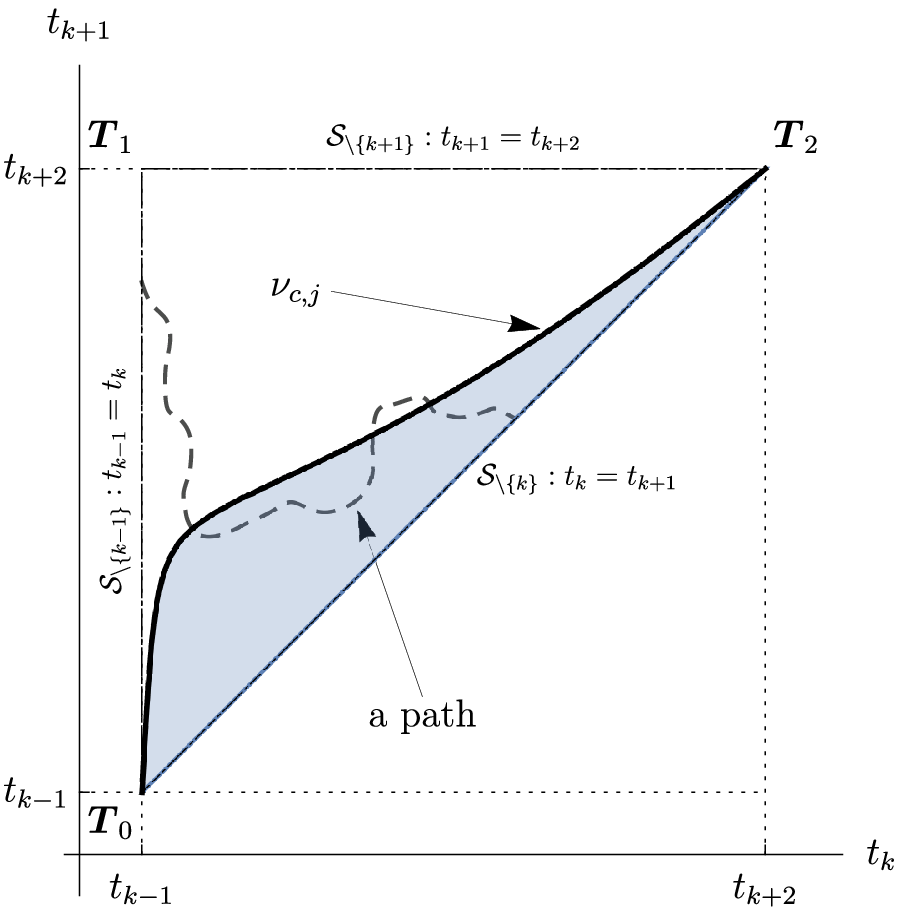}
  \end{minipage}
\caption{The variety $\varc$ separates an open edge on the polygon $\trial\triaul\triaur$ from the  open line segment $\trial\triaur$. }
\label{figLemcKKpOnePf}
\end{figure}

\begin{figure}[htb]
	\centering
  \begin{minipage}{0.47\textwidth}
    \centering \includegraphics[width=\textwidth]{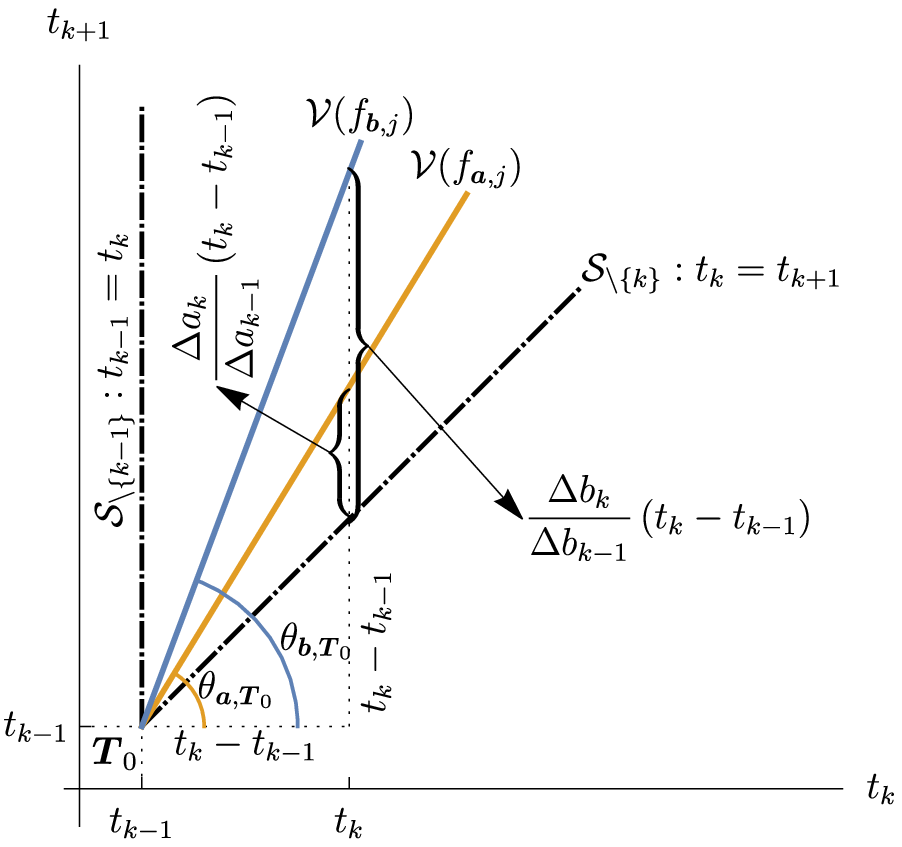}
  \end{minipage}
  \begin{minipage}{0.47\textwidth}
    \centering \includegraphics[width=\textwidth]{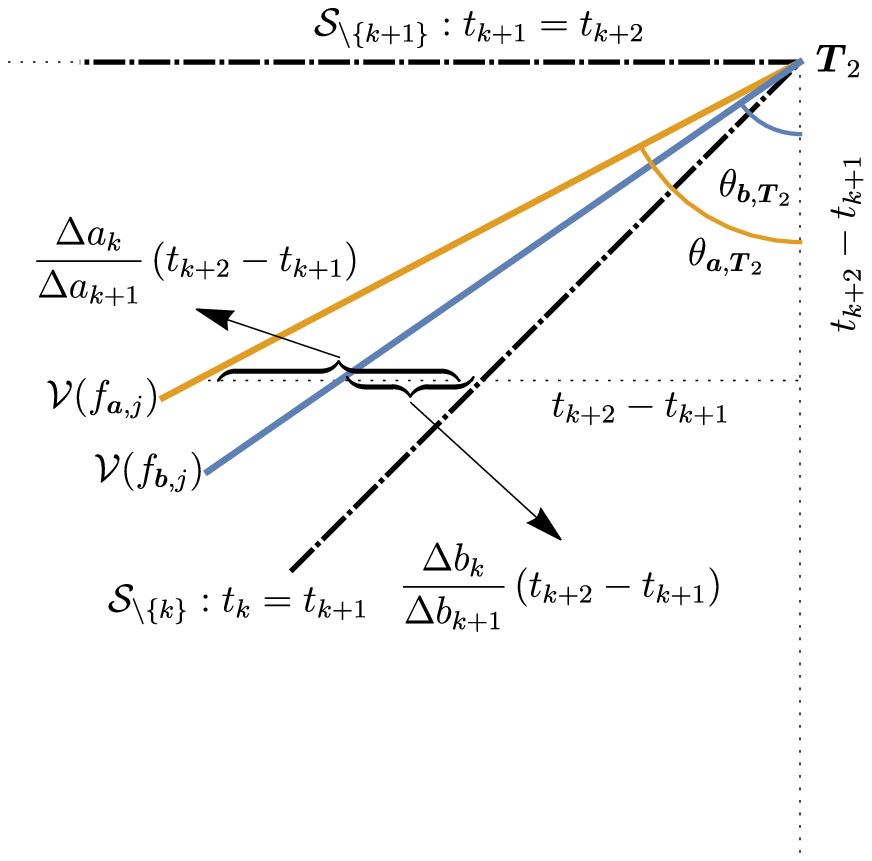}
  \end{minipage}
\caption{Interlacing of varieties $\vara$, and $\varb$ at the vertices $\trial$ and  $\triaur$. }
\label{figLemabKKpOnePf}
\end{figure}

\begin{figure}[htb]
	\centering
 \begin{minipage}{0.45\textwidth}
    \centering \includegraphics[width=\textwidth]{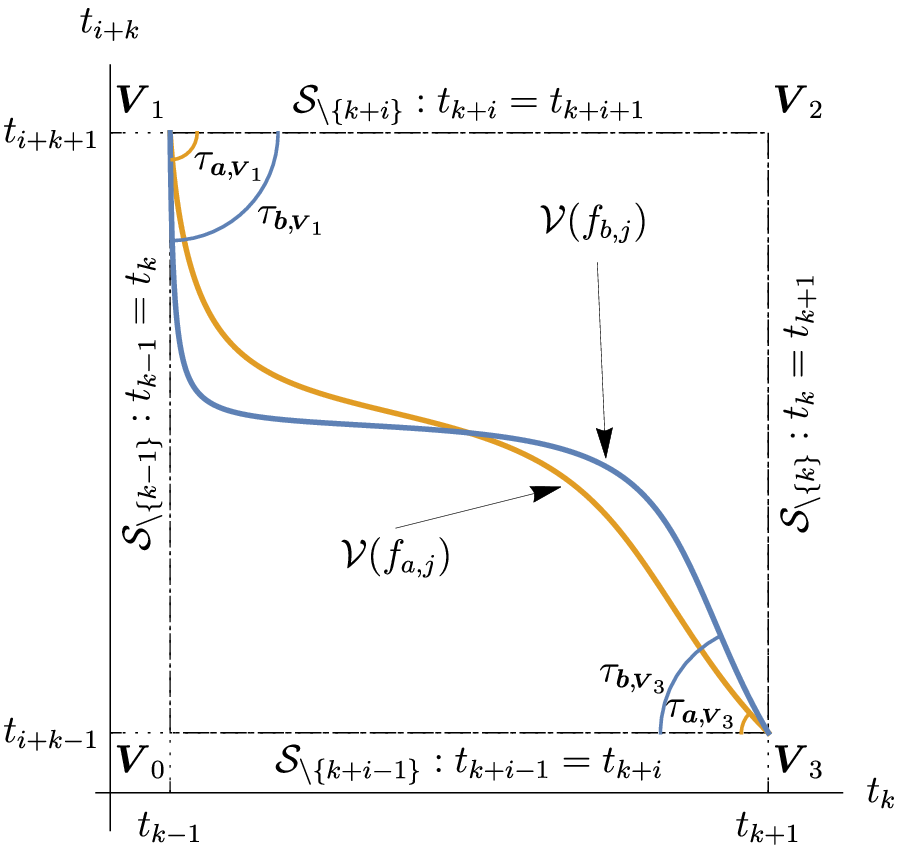}
 \end{minipage}
 \begin{minipage}{0.45\textwidth}
   \centering \includegraphics[width=\textwidth]{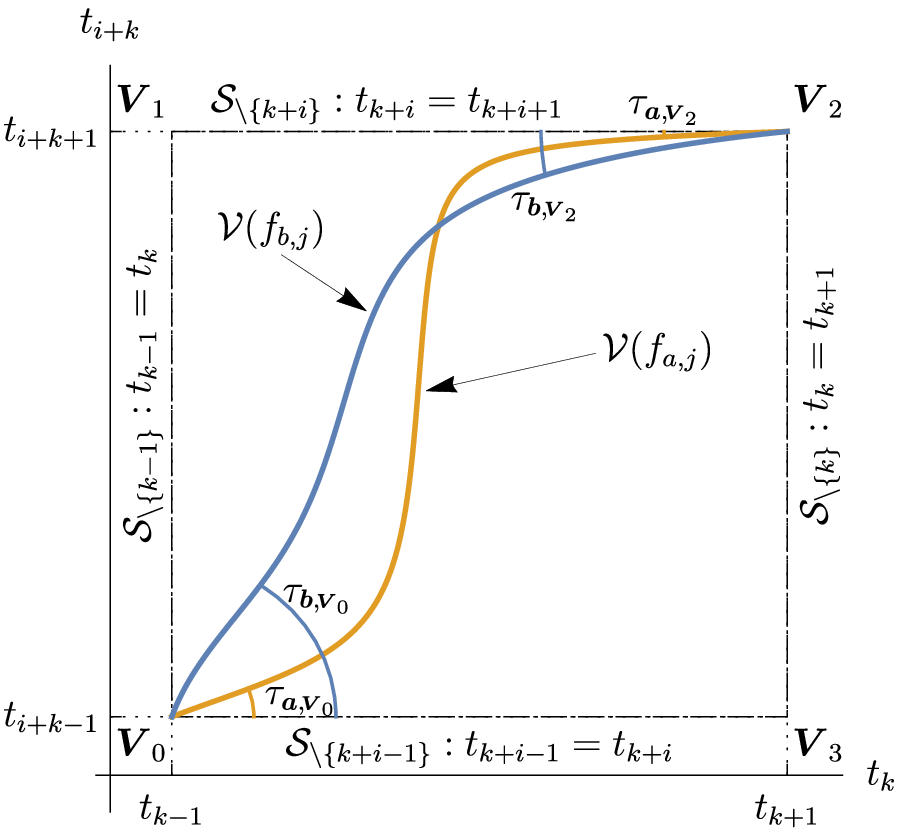}
  \end{minipage}
	\caption{Interlacing of the varieties $\vara$ and $\varb$ in $\interior{\qrial\qriaul\qriaur\qriar}$ for even and odd $i$, left and right respectively.}
\label{figLemabKKpiPf}
\end{figure}

With $j$ fixed, let us choose $k - 1, \, \ j \le k \le n+j$,  as a candidate for $\lstarp{\vv{c},j}$. Following  \lemnref{lem:lemma1} we have to analyze $\fc$ near the faces
\begin{equation}\label{eqnigh1}
\underbrace{\smfp{j-1,k-1}, \dots,\smfp{k-3,k-1},\smfp{k-2,k-1}}_{k-j},
\end{equation}
and
\begin{equation}\label{eqnigh2}
\underbrace{\smfp{k-1,k},\smfp{k-1,k+1},\dots,\smfp{k-1,n+j-1}}_{n+j-k}
\end{equation}
since $\fc$   vanishes at these $(2n-4)$-simplices. 
Let us consider now \eqref{eqnigh2} only since the inspection of  \eqref{eqnigh1} follows by  symmetry. 
The faces will be analyzed  with help of an intersection of $\sm$ and the ${\left(t_{k},t_{k+i}\right)}$-plane, $1 \le i \le n + j - 1 - k$, 
with the rest of the components $\vv{t} \in \interior \sm$ considered as parameters.
Two different  generic cases are to be examined only, $i=1$ yields the intersection as a triangle
$$
\left\{ \left(  t_k, t_{k+1} \right)^T \in \RR^2 | \, t_{k-1} \le t_k \le t_{k+1} \le t_{k+2} \right\},
$$
and the rest of $i$ as a rectangle
$$
\left\{ \left(  t_k, t_{k+i}\right)^T \in \RR^2|\, t_{k-1}  \le t_k \le t_{k+1}, \, t_{k+i-1} \le t_{k+i} \le  t_{k+i+1} \right\}.
$$
As the first one, we take a look at the intersection of $\sm$ and the $(t_k, t_{k+1})$-plane. This intersection  is a triangle, determined by the vertices
\begin{equation}\label{tdefs}
\trial := \vdva{ t_{k-1}}{t_{k-1}}, \  \triaul := \vdva{t_{k-1}}{t_{k+2}},\  \triaur := \vdva{t_{k+2}}{t_{k+2}},
\end{equation}
\ie particular points in 
$$
\smfp{k-1,k}, \smfp{k,k+1}, \smfp{k-1,k+1},
$$
and edges as particular line segments in
$$
\smfp{k-1}, \smfp{k}, \smfp{k+1}.
$$
If $\vv{t} \in \interior \sm$ moves, the functions $\fc$ as well as the actual positioning of the points \eqref{tdefs} in $\RR^{2n-2}$ change  but the crucial relations that will be proved in the following lemmas remain intact.
\begin{lem}\label{lemkkp1} 
Suppose that $j \le k \le n + j  - 2$, and
\begin{equation}\label{asslemkkp1} 
\Delta c_{k-1} > 0,  \ \Delta c_{k} > 0 , \ \Delta c_{k+1} > 0.
\end{equation}
Suppose that the components of \, $\vv{t} \in \sm$ satisfy
\begin{equation}\label{pogvtdif}
t_{j-1} < t_j < \dots < t_{k-1} < t_{k+2} < t_{k+3} < \dots < t_{n+j}.
\end{equation}
At the triangle $\trial\triaul\triaur$, the variety $\varc$ has the following properties \figref{figLemcKKpOnePf}:
\begin{enumerate} [label=(\roman*)]
\item near the vertices  $\triac{\ell}$ it could be expressed as a continuous  function of $t_{k}$ (or $t_{k+1}$).  At $\triaul$ it only touches the triangle from the outside, and at  $\trial$ as well as at $\triaur$ it continues in the triangles interior,
\item any path connecting the open line segment $\trial\triaur$ with the open edge on the polygon $\trial\triaul\triaur$ crosses the variety odd times,
\item it lies entirely inside the triangle except for the vertices $\tria_i$, and it continuously connects $\trial$ and $\triaur$.
\end{enumerate}
\end{lem}
\begin{pf}
The assumptions \eqref{asslemkkp1}, \eqref{pogvtdif},  and \lemnref{lem:lemma1} imply that the function $\fc$ is of signs $(-1)^{n+j-k}$, $(-1)^{n+j-k}$, and $(-1)^{n+j-k+1}$ at the open line segments $\trial\triaul$, $\triaul\triaur$, and $\triaur\trial$ respectfully.  The expansion \eqref{fcexp} shows that we may apply the implicit function theorem to $\fc$ to obtain $\varc$ at $\triac{\ell}$  locally as a continuous function of $t_{k}$ or $t_{k+1}$.
At $\triaur$, the significant part of the expansion reads
\begin{align*}
& \linp{k-1,k,k+2,k+1} = \\
& \quad = - \omega '_{k,k+1}\left(t_{k-1}\right) \Delta c_{k+1} \left(t_k-t_{k-1}\right)
+ \omega '_{k,k+1}\left(t_{k+2}\right)  \Delta c_{k-1} \left(t_{k+2}-t_{k+1}\right) .
\end{align*}
Since
$
\sign \ \omega '_{k,k+1}\left(t_{k+2}\right) = - \sign \  \omega '_{k,k+1}\left(t_{k-1}\right) = (-1)^{n+j-k},
$
the first assertion is verified at $\triaur$, and similarly at $\trial$. At $\triaul$, the first order expansion can vanish only if
$\sign \Delta t_{k-1} = - \sign \Delta t_{k+1} $, and the touch is verified too.
Further, the variety $\varc$ must block any path from the open line segment $\trial\triaur$ to $\triaul$, so it must continuously connect $\trial$ and $\triaur$, and the rest of the  assertions is confirmed.
\qed
\end{pf}
In the second case, we consider the rectangular intersection of $\sm$ and${ (t_k, t_{k+i})\textrm{-plane}}$ \  where $2 \le i$, $k+i \le n+j - 1$. Its edges correspond to particular line segments of the faces
$$
\smfp{k-1}, \smfp{k+i}, \smfp{k},  \smfp{k+i-1},
$$
and its vertices
$$
\qrial := \vdva{ t_{k-1}}{t_{k+i-1}}, \  \qriaul := \vdva{t_{k-1}}{t_{k+i+1}},\  \qriaur := \vdva{t_{k+1}}{t_{k+i+1}}, \ \qriar := \vdva{ t_{k+1}}{t_{k+i-1}},
$$
are particular points of their intersections
$$
 \smfp{k-1,k+i-1}, \smfp{k-1,k+i}, \smfp{k,k+i}, \smfp{k,k+i-1}.
$$
We have the following observation, very similar to \lemnref{lemkkp1}.
\begin{lem}\label{lemkkpi} 
Suppose that $j \le k$, \ $k+1 < k + i \le n+j- 1$, and
\begin{equation}\label{asslemkkpi}
\Delta c_{k-1} > 0,  \ \Delta c_{k} > 0 , \ \Delta c_{k+i-1} > 0, \ \Delta c_{k+i} > 0.
\end{equation}
Suppose that the components of  $\vv{t}  \in \sm$ satisfy
\begin{equation}\label{pogvtdifa}
t_{j-1} <  \dots < t_{k-1} < \ t_{k+1} <\dots < t_{k+i-1} < t_{k+i+1} < \dots <  t_{n+j}.
\end{equation}
At the rectangle $\qrial\qriaul\qriaur\qriar$ the variety $\varc$ has the following properties \figref{figLemKKpiPf}. Near the vertices $\qriac{\ell}$ it  could be expressed as a continuous function of $t_{k}$ (or $t_{k+i}$), and for an even $i$ one has
\begin{figure}[htb]
	\centering
  \begin{minipage}{0.45\textwidth}
    \centering \includegraphics[width=\textwidth]{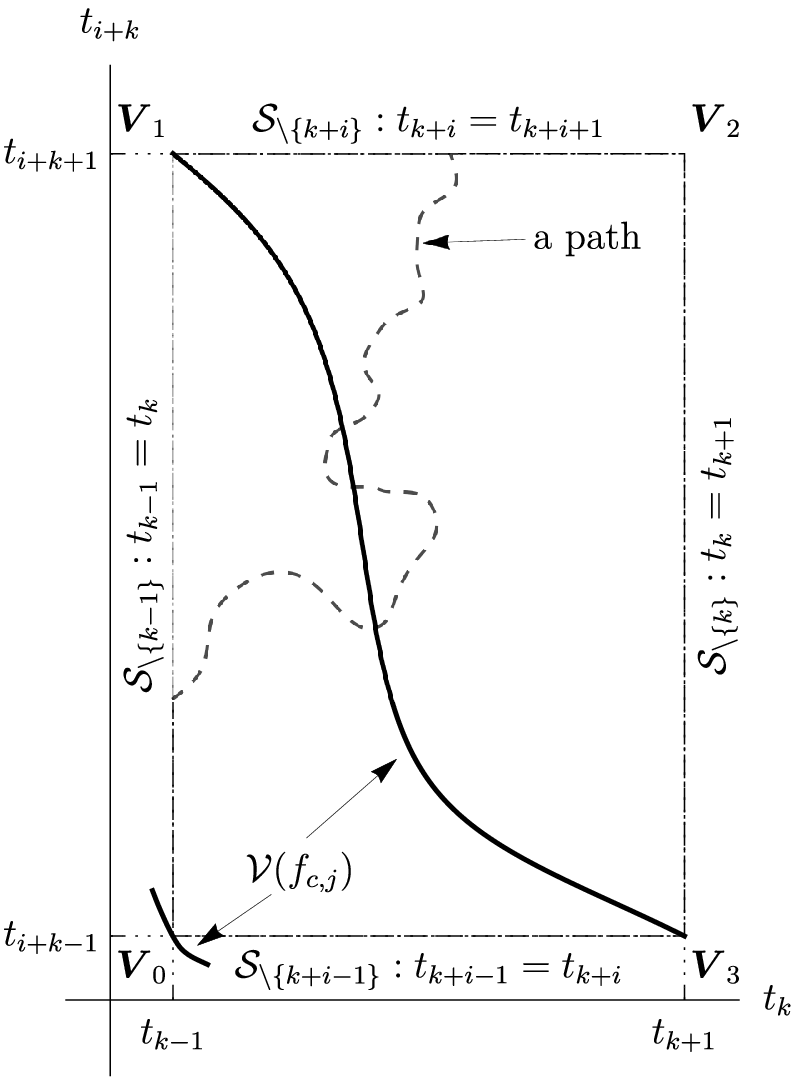}
  \end{minipage}
  \begin{minipage}{0.45\textwidth}
    \centering \includegraphics[width=\textwidth]{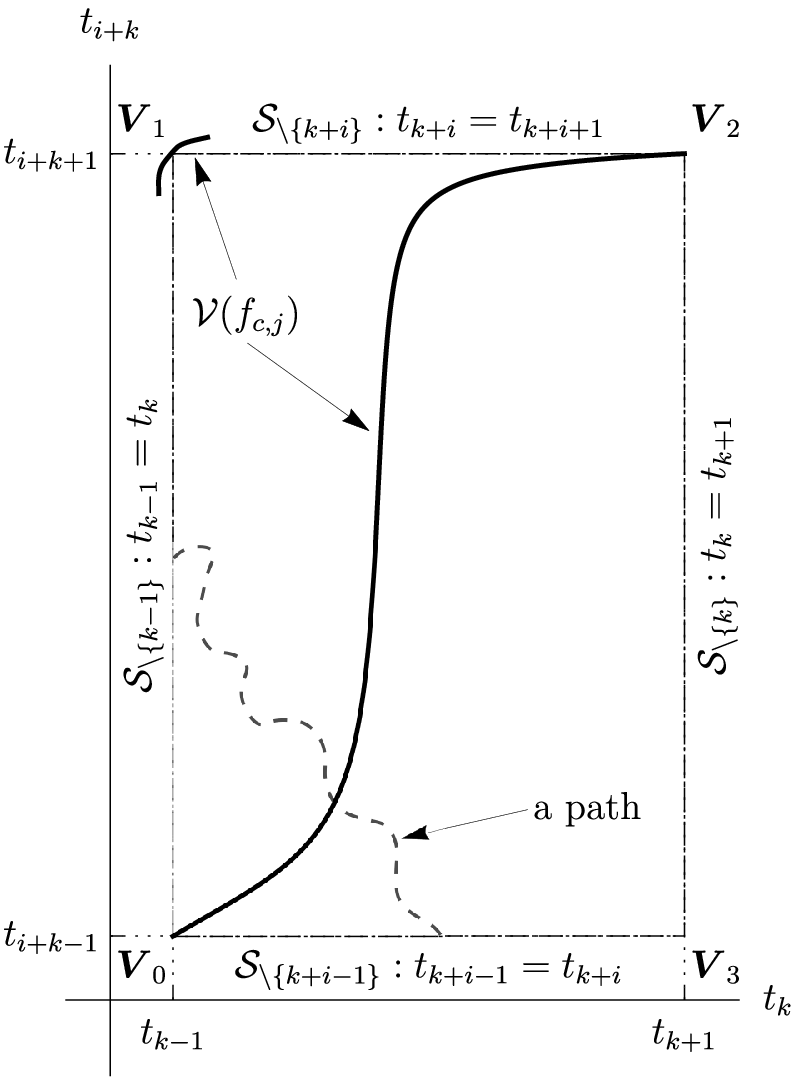}
  \end{minipage}
 \caption{The variety $\varc$ separates an open edge on the  polygon $\qriar\qrial\qriaul$ from an open edge on $\qriaul\qriaur\qriar$ (left, $i$ even), and an open edge on $\qrial\qriaul\qriaur$ from an open edge on $\qriaur\qriar\qrial$ (right, $i$ odd). }
\label{figLemKKpiPf}
\end{figure}
\begin{enumerate} [label=(\roman*)]
\item at $\qrial$ and $\qriaur$ it only touches the rectangle, but at  $\qriaul$ and $\qriar$ it continues to its interior,
\item any path connecting an open edge on the  polygon $\qriar\qrial\qriaul$ to an open edge on $\qriaul\qriaur\qriar$ crosses the variety an odd number of times,
\item $\varc$ lies entirely inside the rectangle except for the vertices, and it continuously connects $\qriaul$ and $\qriar$.
\end{enumerate}
If $i$ is odd, the role of $\qrial$ and $\qriaul$, as well as $\qriaur$ and $\qriar$, is reversed.
\end{lem}
\begin{pf}
The assumption \eqref{pogvtdifa} allows us to use the expansion \eqref{lindef} again, and we obtain significant expansion parts at four rectangle corners as,
\begin{align*}
& \linp{k-1,k,k+i-1,k+i} = 
 (-1)^{i} \left(-  \text{$\Delta $c}_{k+i-1} \omega '_{k,k+i}\left(t_{k-1}\right) \left(t_{k}-t_{k-1}\right)\right. + \nonumber \\
& \quad\quad + \left. \text{$\Delta $c}_{k-1}
  \omega '_{k,k+i}\left(t_{k+i-1}\right)  \left(t_{k+i}-t_{k+i-1}\right)  \right), \quad  \textrm{at} \ \qrial, 
\end{align*}
\begin{align*}
& \linp{k-1,k,k+i+1,k+i} = 
(-1)^i \left(- \text{$\Delta $c}_{k-1}  \omega '_{k,k+i}\left(t_{k+i+1}\right)\left(t_{k+i+1}-t_{k+i}\right)\right. + \nonumber \\
& \quad\quad + \left.  \text{$\Delta   $c}_{k+i} \omega '_{k,k+i}\left(t_{k-1}\right) \left(t_{k}-t_{k-1}\right) \right), \quad  \textrm{at} \ \qriaul, 
\end{align*}
\begin{align*}
& \linp{k+1,k,k+i+1,k+i} = 
 (-1)^i \left( -  \text{$\Delta $c}_{k+i} \omega '_{k,k+i}\left(t_{k+1}\right) \left(t_{k+1}-t_{k}\right) \right. + \nonumber \\
& \quad\quad + \left. \text{$\Delta $c}_k
   \omega '_{k,k+i}\left(t_{k+i+1}\right)  \left(t_{k+i+1}-t_{k+i}\right) \right), \quad  \textrm{at} \ \qriaur, 
\end{align*}
\begin{align*}
& \linp{k+1,k,k+i-1,k+i} = 
(-1)^i \left( \text{$\Delta $c}_{k+i-1} \omega '_{k,k+i}\left(t_{k+1}\right) \left(t_{k+1}-t_{k}\right) \right. - \nonumber \\
& \quad\quad - \left. \text{$\Delta $c}_k
   \omega '_{k,k+i}\left(t_{k+i-1}\right) \left(t_{k+i}-t_{k+i-1}\right)  \right), \quad  \textrm{at} \ \qriar. 
\end{align*}
Note also
\begin{align}\label{signomkpi}
& \sign \ \omdcp{k,k+i}{t_{k+1}} =  - \sign \ \omdcp{k,k+i}{t_{k-1}} = (-1)^{n+j-k}, \nonumber \\[-6pt]
& \\[-6pt]
& \sign \ \omdcp{k,k+i}{t_{k+i-1}} =  - \sign \ \omdcp{k,k+i}{t_{k+i+ 1}} = (-1)^{n+j-k+i}. \nonumber
\end{align}
With the assumption \eqref {asslemkkpi} it is now straightforward to verify the assertions by arguments used already in the proof of \lemnref{lemkkp1}. \qed
\end{pf}
To investigate the system \eqref{system-pol-form} further, we have to make use of the  distinction between $\fa$, and $\fb$. Recall \lemnref{lemkkp1} and \fignref{figLemcKKpOnePf}. The variety $\varc$ continuously connects $\trial$ and $\triaur$. Let us denote by $\thtcnic$ the angle between abscissa direction and $\varc$ at $\trial$, and let $\thtcdva$ be the angle between $\varc$ and ordinate direction at $\triaur$.
%
%
%
\begin{lem}\label{lemkkp1ab} 
Suppose assumptions of \lemnref{lemkkp1} are satisfied  and  suppose additionally
\begin{equation}\label{asslemkkp1ab}
\dt{\deap{k-1+\ell}}{\deap{k+\ell}}{\debp{k-1+\ell}}{\debp{k+\ell}}  > 0, \quad \ell = 0,1.
\end{equation}
Then \figref{figLemabKKpOnePf}
\begin{align}\label{inetan}
& 1 < \tan{\thtanic} = 1 + \frac{\deap{k}}{\Delta a_{k-1}} < \tan{\thtbnic} = 1 + \frac{\debp{k}}{\Delta b_{k-1}} < \infty, \nonumber\\[-6pt]
& \\[-6pt] 
&1 < \tan{\thtbdva} = 1 + \frac{\debp{k}}{\Delta b_{k+1}} < \tan{\thtadva} = 1 + \frac{\deap{k}}{\Delta a_{k+1}} < \infty. \nonumber
\end{align}
The number of intersections of the varieties $\vara$ and $\varb$ in $\interior{\trial\triaul\triaur}$ is odd.
\end{lem}
\begin{pf}
At $\trial$, the variety $\varc$ is determined from the significant part of the $\fc$ expansion \eqref{fcexp} 
\begin{align*}
& \linp{k-1,k,k-1,k+1} = 0 = \nonumber\\
& \quad\quad =    \omega '_{k,k+1}(t_{k-1}) \left( (\Delta c_{k-1} + \Delta c_{k}) \left(t_k - t_{k-1}\right)  -  \Delta c_{k-1} \left(t_{k+1} - t_{k-1}\right) \right),
\end{align*}
so
$$
 \tan{\thtcnic} = 1 + \frac{\Delta c_{k}}{\Delta c_{k-1}},
$$
 and similarly
$$
 \tan{\thtcdva} = 1 + \frac{\Delta c_{k}}{\Delta c_{k+1}}
$$
at $\triaur$.  The inequalities \eqref{inetan} follow then from \eqref{asslemkkp1} and \eqref{asslemkkp1ab}. So varieties $\vara$ and $\varb$  interlace, 
and the proof is completed.  \qed
\end{pf}
%
%
\begin{figure}[htb]
	\centering
 \begin{minipage}{0.45\textwidth}
    \centering \includegraphics[width=\textwidth]{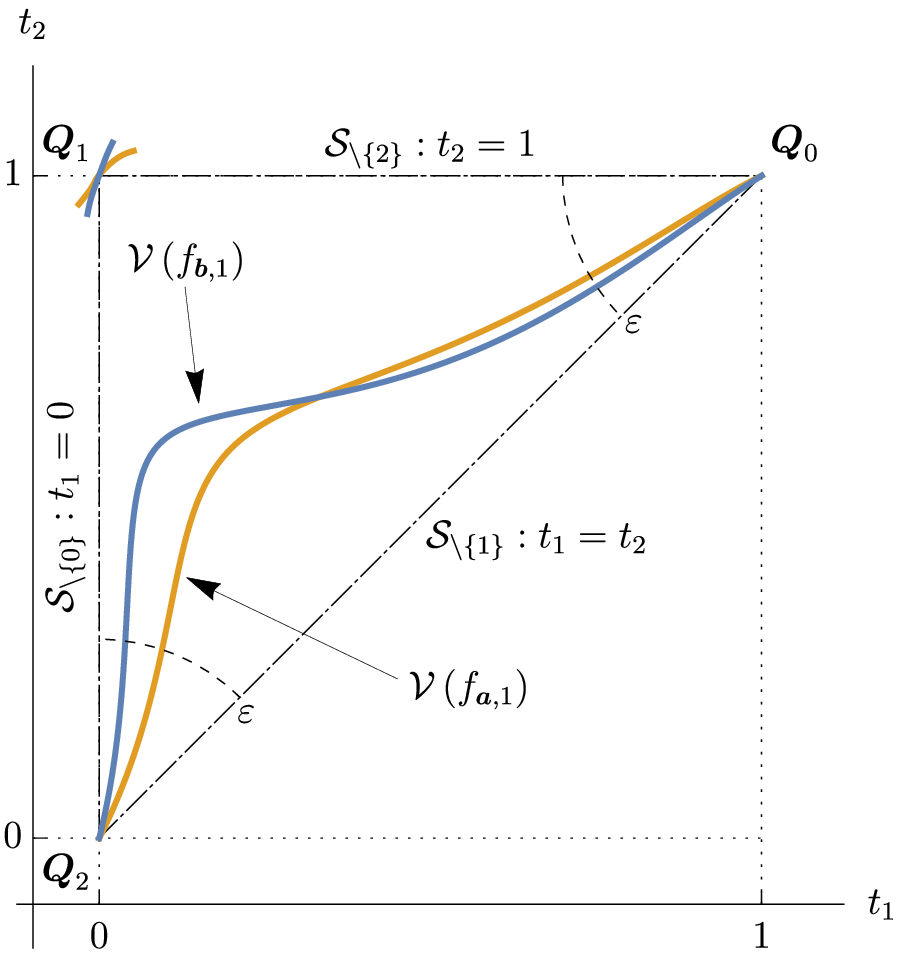}
 \end{minipage}
\begin{minipage}{0.45\textwidth}
    \centering \includegraphics[width=\textwidth]{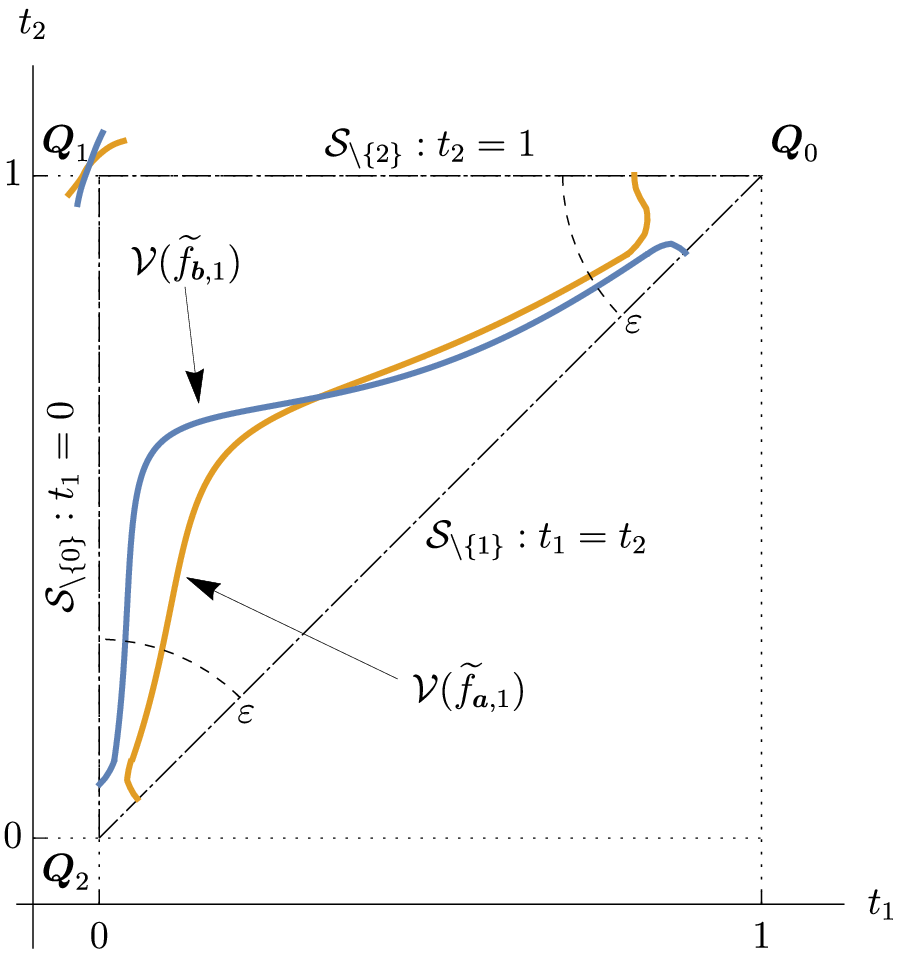}
 \end{minipage}
\caption{The case $n=2$ before and after modifications of $\fa \to \tfa$ and $\fb \to \tfa$ based upon  \lemnref{lemkkp1ab} and \lemnref{lemmodstep} (the image left and right respectively).}
 \label{figRemNisTwo}
\end{figure}
\begin{rmk}\label{rem:remark4}
In the particular case, \ie $n=2$, the proof of \thmnref{theorem1} is completed. 
If the assumptions are satisfied, \lemnref{lemkkp1ab} in hand with the implicit function theorem reveals  that modifications ${\funpp{\vv{a}}{1} \to  \tfunpp{\vv{a}}{1}}$, ${\funpp{\vv{b}}{1} \to  \tfunpp{\vv{b}}{1}}$  which satisfy \lemnref{lemmodstep} can be carried over at $\qpp{2}$ and $\qpp{0}$ in some $\varepsilon$ neighborhood \figref{figRemNisTwo}.  The 
function $\tfunpp{\vv{a}}{1}$  is of opposite sign at $\smfp{0}$ and $\qpp{0}$,  so is $\tfunpp{\vv{b}}{1}$ at $\smfp{2}$ and $\qpp{0}$.
Thus 
$$ 
{\lstarp{{\vv{a},1}} = 0}, \  {\lstarp{{\vv{b},1}} = 2}, \ \  \textrm{and} \  {\ldstar=1}.
$$
After removing the isolated zero at $\qpp{1}$ we have
$$
\sign \tfunpp{\vv{a}}{1}\left(\qpp{1}\right)  = \sign \tfunpp{\vv{b}}{1}\left(\qpp{1}\right) = 1,
$$ 
but nowhere at $\smfp{1}$ both functions are positive simultaneously. So the existence of the solution follows from \citep[ Theorem 2.1]{Vrahatis-2016}. Of course, there are much shorter ways to handle $n = 2$, even the closed form solution is available,
\begin{align}\label{eq:nI2clf}
& \dPijS{i}{k} :=\dtP{\dTp{i}}{\!\!\!\dTp{k}} >0,\quad i=0,1; k =1,2,\\
& t_1 = t_0 + \frac{ \dPijS{0}{1} \dPijS{0}{2}}{\dPijS{0}{1} \left(\dPijS{1}{2}+\dPijS{0}{2}\right)
+ \sqrt{ \dPijS{0}{1} \dPijS{1}{2} \left(\dPijS{0}{1}+\dPijS{0}{2}\right)  \left( \dPijS{1}{2} + \dPijS{0}{2} \right)}}, \nonumber \\
& t_2 = t_0 + \frac{\dPijS{0}{1} \dPijS{1}{2} + \sqrt{\dPijS{0}{1} \dPijS{1}{2} \left(\dPijS{0}{1} + \dPijS{0}{2}\right)  \left( \dPijS{1}{2} + \dPijS{0}{2} \right)}}{ \dPijS{1}{2} \left(\dPijS{0}{1}+\dPijS{0}{2}\right) +\sqrt{ \dPijS{0}{1} \dPijS{1}{2} \left(\dPijS{0}{1}+\dPijS{0}{2}\right) \left( \dPijS{1}{2} + \dPijS{0}{2} \right)}}. \nonumber
\end{align}
Nevertheless, this remark suggests how to complete the general case too.
\end{rmk}
Since the case $n = 2$ is covered, we shall assume from now on $n \ge 3$. We return now to \lemnref{lemkkpi}, with sharpened assumptions. Let us denote by $\ctau{\bfm{c}}{\ell}$ the angle between abscissa direction and $\varc$ at $\qria_\ell$.
\begin{lem}\label{lemkkpiab} 
Suppose assumptions of \lemnref{lemkkpi} are satisfied and  suppose additionally
\begin{equation}\label{poglem6}
\dt{\deap{k-1+\ell}}{\deap{k+i+\ell}}{\debp{k-1+\ell}}{\debp{k+i+\ell}} > 0,
\quad \ell = 0,1.
\end{equation}
Then the angles $\tauc{\qria_\ell}$ at vertices $\qria_\ell$ satisfy  \figref{figLemabKKpiPf}
\begin{equation*}
0 < \tan{\tauaena}  < \tan{\taubena} < \infty, \quad 0 < \tan{\tauatri}  < \tan{\taubtri}  < \infty,
\end{equation*}
if $i$ is even, and
\begin{equation*}
0 < \tan{\tauanic}  < \tan{\taubnic} < \infty, \quad 0 < \tan{\tauadva}  < \tan{\taubdva}  < \infty,
\end{equation*}
for $i$ that is odd. The number of intersections of the varieties $\vara$ and $\varb$ in $\interior{\qrial\qriaul\qriaur\qriar}$ is odd.
\end{lem}
\begin{pf}
If $i$ is even, the linearized equation
$$
\linp{k-1,k,k+i+1,k+i} = 0
$$
at $\qriaul$ yields
\begin{equation*}
\tan \taucena =\frac{ \omega '_{k,k+i}\left(t_{k-1}\right)}{ \omega '_{k,k+i}\left(t_{k+i+1}\right)} \frac{\Delta c_{k+i}}{\Delta c_{k-1}} ,
\end{equation*}
and 
$$
\linp{k+1,k,k+i-1,k+i} = 0
$$
similarly
\begin{equation*}
\tan \tauctri = \frac{ \omega '_{k,k+i}\left(t_{k+1}\right)}{ \omega '_{k,k+i}\left(t_{k+i-1}\right)} \frac{\Delta c_{k+i-1} }{\Delta c_k}
\end{equation*}
at $\qriar$.
The asserted inequalities follow then from \eqref{asslemkkpi}, \eqref{signomkpi},  and \eqref{poglem6}. We omit the odd case $i$ part of the proof. 
\qed
\end{pf}
Let us  verify now that \eqref{eqthmsupx} and \eqref{eqthmsupz}  imply the determinant sign assumptions of  the previous lemmas of this section. 
\begin{lem}\label{lemma6a}
Let $k,\ell,i$ be indices that satisfy $0\le k < \ell < i < 2n-1$. Suppose that
$$
	\deap{k} >0 , \debp{k} >0, \deap{\ell} >0, \debp{\ell} >0, \deap{i} >0,
$$
and
$$
\dt{\deap{k}}{\deap{\ell}}{\debp{k}}{\debp{\ell}} > 0, \
 \dt{\deap{\ell}}{\deap{i}}{\debp{\ell}}{\debp{i}} > 0
$$
Then
\begin{equation}\label{detlipos}
\dt{\deap{k}}{\deap{i}}{\debp{k}}{\debp{i}} > 0.
\end{equation}
\end{lem}
\begin{pf}
If the assumptions  hold, but not \eqref{detlipos}, we observe
$$
\debp{k} \ge \debp{i} \frac{\deap{k}}{\deap{i}} > \debp{\ell}  \frac{\deap{i}}{\deap{\ell}} \frac{\deap{k}}{\deap{i}} >  \debp{k}  \frac{\deap{\ell}}{\deap{k}} \frac{\deap{k}}{\deap{\ell}} =  \debp{k},
$$
a contradiction. \qed
\end{pf}
Let the assumptions of \thmnref{theorem1} be satisfied. A brief look reveals that 
$$
	\deap{k}  \Delta b_{i} - \Delta a_{i} \debp{k}> 0, \quad  0 \le k < i \le 2n-1,
$$
clearly follows from  \lemnref{lemma6a}, so all the assumptions in \lemnref{lemkkp1}, \lemnref{lemkkpi}, \lemnref{lemkkp1ab}, and   \lemnref{lemkkpiab} are met too. 
\begin{figure}[htb]
	\centering
 \begin{minipage}{0.45\textwidth}
    \centering \includegraphics[width=\textwidth]{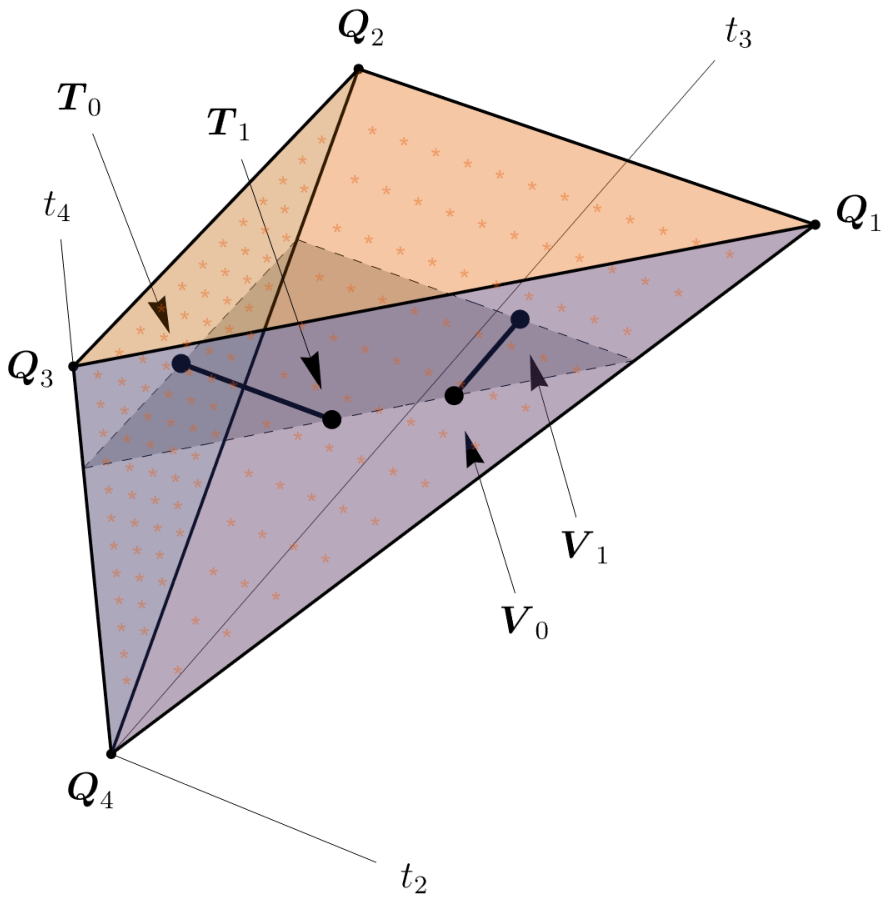}
 \end{minipage}
 \begin{minipage}{0.45\textwidth}
   \centering \includegraphics[width=\textwidth]{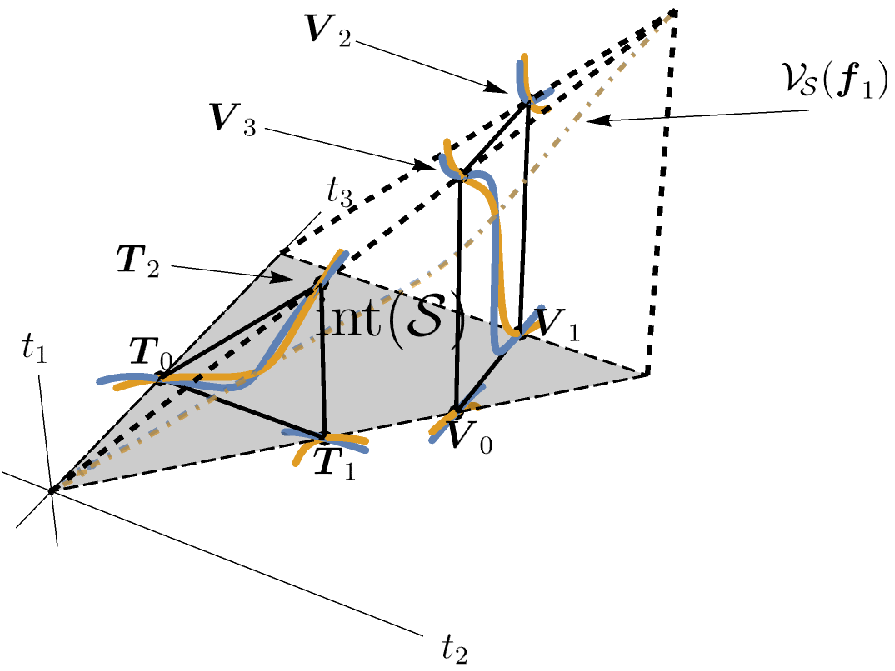}
  \end{minipage}
	\caption{The cubic case, $\sm \in \RR^4$, and the functions $\fFunj{1} = \left( \faap{1},\fbbp{1}\right)^T$ considered. The boundary simplex $\smfp{0}$ in the
	$(t_2,t_3,t_4)^T$-subspace, with a particular choice of $\trial\triaul$ and $\qrial\qriaul$ based upon the same choice of $t_4$ (left), and the corresponding intersections discussed in \lemnref{lemkkp1ab} or \lemnref{lemkkpiab} 
in $(t_1,t_2,t_3)^T$-subspace of $\interior{\sm}$ (right). The dash-dotted curve shows $\varisp{\fFunj{1}}$   
for the particular choice of~$t_4$. It connects $\qpp{3}\qpp{4}$ and $\qpp{0}\qpp{4}$.}
\label{fignIsthreeL}
\end{figure}
Let us introduce further discussion by a cubic example \figref{fignIsthreeL}.  In this case, $\sm$ lives in $\RR^4$, and only projections are available.
The left image shows the boundary simplex $\smfp{0}$ in $(t_2,t_3,t_4)^T$-subspace. The line segments  $\trial\triaul$ and $\qrial\qriaul$, referenced in \lemnref{lemkkp1}, \lemnref{lemkkp1ab}, and \lemnref{lemkkpi},   \lemnref{lemkkpiab} respectively, are  obtained at $t_4 = \const$. The functions $\faap{1}$, 
$\fbbp{1}$ are negative inside $\smfp{0}$, and they change the sign when crossing through the faces $\smfp{0,1}\ (  = \cop{\qpp{2}\qpp{3}\qpp{4} })$ and $\smfp{0,3}$ to the simplices $\smfp{1}$ and $\smfp{3}$ respectively. The right image shows the corresponding $\sm$ projection at $t_4 = \const$. 
Curves shown  at $\trial\triaul\triaur$ or $\qrial\qriaul\qriaur\qriar$ are parts of $\varp{\faap{1}}$ and $\varp{\fbbp{1}}$.
Near  $\smfp{0,1}$ or $\smfp{0,1}$ they both could  be expressed as continuous functions of variables other that ${t_1}$. The variety  $\varisp{\fFunj{1}}$ as  dash-dotted curve connects points at $\smfp{0,1,2} = \qpp{3}\qpp{4} $ and $\smfp{1,2,3} = \qpp{0}\qpp{4}$. Let us denote
$
\vtj{\ell} := \left(t_1,\dots,t_{\ell-1},t_{\ell+1},\dots,t_{2n-2}\right)^T.
$
\begin{lem}\label{lem:fcatbSj}
Suppose the assumptions of \thmnref{theorem1} are satisfied. 
Let $\mathcal{M}$ stand either for $\smfp{j-1,j-1+k}$ or $\smfp{n+j-1-k,n+j-1}$, with ${1 \le k \le n}$. Then 
\begin{enumerate}[label=(\roman*)]\setlength\itemsep{-4pt}
\item if $k$ is even, $\interior{\mathcal{M}} \cap \varisp{\fc} = \emptyset$,  \label{en:fcatbSji}
\item if $k$ is odd, $\mathcal{M} \subset \varisp{\fc}$. The variety  could be, close enough to $\mathcal{M}$,  expressed as a continuous function
$$
\varisp{\fc}: \RR^{2n-3} \to \RR: \vtj{\ell} \mapsto \varisp{\fc}\left(\vtj{\ell} \right),
$$
with $\ell = j$ for the first case of $\mathcal{M}$, and $\ell = n+j-1$ for the second one. \label{en:fcatbSjii}
\item $\interior{\mathcal{M}} \cap \varisp{\fFunj{j}} = \emptyset$, \label{en:fcatbSjiii}
\item $\smfp{j-1} \cap \varisp{\fFunj{j}} = \smfp{j-1,j,j+1}$, 
$\smfp{n+j-1} \cap \varisp{\fFunj{j}} = \smfp{n+j-3,n+j-2,n+j-1}$, and  $\varisp{\fFunj{j}}$ connects $\smfp{j-1,j,j+1}$ with $\smfp{n+j-3,n+j-2,n+j-1}$. \label{en:fcatbSjiiii}
\end{enumerate} 
\end{lem}
\begin{pf}
\lemnref{lemma6a} allows one to apply previous lemmas of this section. The assertion \ennref{en:fcatbSji} is confirmed by \lemnref{lemkkpi} and \lemnref{lemkkp1ab}. So is \ennref{en:fcatbSjii}, but only on $\interior{\mathcal{M}}$. Since $\fc$ is a polynomial, vanishing at $\boundary{\smfp{j-1}}$ and  $\boundary{\smfp{n+j-1}}$, we may extend the conclusion to the entire $\mathcal{M}$ by continuity. For an even $k$, \ennref{en:fcatbSjiii} follows from \ennref{en:fcatbSjii}, and  \lemnref{lemkkp1ab}, \lemnref{lemkkpiab} imply it for an odd one. From \ennref{en:fcatbSjiii} it follows that $\varisp{\fFunj{j}}$ may include only a part of $\boundary{\mathcal{M}}$,  \ie
parameter vectors $\bfm{t}$ where at least three of the parameters
$$
t_{j-1}, t_{j},\dots,t_{n+j}
$$
 coincide. Here, $t_{j-1}$ and $t_{n+j}$ are assumed to be constant.
Recall  \lemnref{lemkkp1ab}   and \lemnref{lemkkpiab}:  the intersections $\trial\triaul\triaur$ or $\qrial\qriaul\qriaur\qriar$ studied there provided the basic step of the \ennref{en:fcatbSjiii} proof. In particular, $\trial\triaul\triaur$ keeps $\varisp{\fFunj{j}}$ trapped inside the triangle independently of the boundary as long it remains nontrivial, \ie ${t_{j-1} < t_{j+2}}$ for the first case of $\mathcal{M}$, and ${t_{n+j-3} < t_{n+j}}$ for the second one. Thus one may apply continuity again, and the first part of \ennref{en:fcatbSjiiii} is confirmed.
Any intersection considered in \lemnref{lemkkp1ab}  or \lemnref{lemkkpiab} contains an odd number of points of $\varisp{\fFunj{j}}$, so this joint variety must continuously connect ${\smfp{j-1,j,j+1}}$  and $\smfp{n+j-3,n+j-2,n+j-1}$. This concludes the last part of the proof.
\qed
\end{pf} 
%
%
Let us denote
$
\varSBj := \varisp{\fc} \cap \boundary{\smfp{j-1}},
$
and let $\vtZ \in \varSBj$. By \lemnref{lemma6a} there exists  $\coepsZerop{\vtZ} > 0$ such that the variety $\varisp{\fc}$
could be expressed as a continuous function of variables ${\vtj{j}, \vt \in \sm,}$ determined by the implicit function theorem from
$$
\fcp{\ti{j-1},\varisp{\fc},\ti{j+1},\dots,\ti{n+j}} = 0,
$$
as long as they satisfy ${\nrd{\vtZj{j} - \vtj{j}} < \coepsZerop{\vtZ},}$
where $\nrd{\pika}$ denotes the Euclidean norm. 
Since $ \varSBj$ is compact, there exists a smallest, but positive $\coepsZerop{\vtZ}$,
$$
\coepsZeroiS{1} := \min\limits_{\vtZ \in \varSBj}\coepsZerop{\vtZ} > 0,
$$
that holds for both $\fa$ and $\fb$. So the domains of $\varisp{\fa}$ and $\varisp{\fb}$ as  functions  include $\varS{j,\coepsZeroiS{1}}$, with
$$
\varS{j,\coepsZero} := \left\{ \vt \in \sm \, \big| \, \dist{\vtj{j}}{\varSBj}  < \coepsZero \right\}.
$$
Let us elaborate $\varisp{\fFunj{j}}$ near $\smfp{j-1,j,j+1}$ further. Since
\begin{align*}
& \vdmp{t_{j-1}, \dots,t_{k-1},t_{k+1},\dots t_{n+j} } = \\
& = \left\{
\begin{matrix}
\vdmp{t_{j-1}, \dots,t_{k-1},t_{k+1},\dots t_{j+2} } \vdmp{t_{j+3}, \dots, t_{n+j} } \prod\limits_{\mybound{i = j-1}{i \ne k}}^{j+2} \prod\limits_{m=j+3}^{n+j} \left( t_m - t_i \right) , \\[-2pt]
 \quad\quad\quad\quad\quad\quad\quad\quad\quad\quad\quad\ \
 \quad\quad\quad\quad\quad\quad\quad\quad\quad  
 j\le k \le j+2, \nonumber\\
 \\
\vdmp{t_{j-1}, \dots t_{j+2} } \vdmp{t_{j+3}, \dots,t_{k-1},t_{k+1},\dots, t_{n+j} } \prod\limits_{i = j-1}^{j+3} \prod\limits_{\mybound{m=j+4}{m \ne k}}^{n+j} \left( t_m-t_i \right),\\ 
 \quad\quad\quad\quad\quad\quad\quad\quad\quad\quad\quad\ \
 \quad\quad\quad\quad\quad\quad\quad\quad\quad  
 j+3\le k \le n+j,
\end{matrix}
\right.
\end{align*}
and $\vdmp{t_{j-1}, \dots t_{j+2} }$ is by the order~$3$ smaller than $\vdmp{t_{j-1}, \dots,t_{k-1},t_{k+1},\dots t_{j+2} } $ there,  $\fc$ in \eqref{fdifdef} simplifies to 
$$
\fc = \left( c_j - c_{j-1}, c_{j+1} - c_{j-1} , c_{j+2} - c_{j-1} \right) \assvecpj + \dots,
$$
where dots denote higher order terms, and
$$
\assvecpj := \big( \vdmp{t_{j-1},  t_{j+1}, t_{j+2} }, - \vdmp{t_{j-1},  t_{j}, t_{j+2} }, \vdmp{t_{j-1},  t_{j}, t_{j+1} }\big)^T .
$$
Thus the main part of the variety $\varp{\fFunj{j}}$ is  determined from 
\begin{align}\label{eq:vass}
 \asstranspj \assvecpj = \bfm{0} , 
\quad \asstranspj := 
\begin{bmatrix}
a_j - a_{j-1}& a_{j+1} - a_{j-1} & a_{j+2} - a_{j-1} \cr
b_j - b_{j-1}& b_{j+1}- b_{j-1} & b_{j+2}- b_{j-1}
\end{bmatrix}. 
\end{align}
Note that \eqref{eq:vass} is just the quadratic case considered in \remnref {rem:remark4}, with indices $0,1,\dots$ been replaced by $j-1,j,\dots$. The equation \eqref{eq:vass} can be written as ${\assvecpj \in \ker\asstranspj}$, where $\ker\asstranspj$ is spanned by a cofactor vector 
\begin{align}\label{eq:bplusdef}
\left(
\begin{matrix}
\ \ \dt{a_{j+1}-a_{j-1}}{a_{j+2}-a_{j-1}}{b_{j+1}-b_{j-1}}{b_{j+2}-b_{j-1}}\\[+10pt]
- \dt{a_{j}-a_{j-1}}{a_{j+2}-a_{j-1}}{b_{j}-b_{j-1}}{b_{j+2}-b_{j-1}}\\[+10pt]
\dt{a_{j}-a_{j-1}}{a_{j+1}-a_{j-1}}{b_{j}-b_{j-1}}{b_{j+1}-b_{j-1}}
\end{matrix}
 \right)
 .
\end{align}
All the determinants in \eqref{eq:bplusdef} are positive by \lemnref{lemma6a}, and the corresponding parameters 
$$
t_{j-1+k} = t_{j-1} + \const_{j,k} \, \xi, \ \const_{j,k} > 0, \quad k= 1,2,3, \ \xi \ge 0, 
$$
could be derived similarly to \eqref{eq:nI2clf}. They depend on the data involved in \eqref{eq:bplusdef} only. Thus we have established 
\begin{align}\label{eq:varabass}
  \varisp{\fFunj{j}}   = & \ \bfm{t}^0 + \xi \, \bfm{h}_j  + \dots, \quad \bfm{t}^0 \in {\smfp{j-1,j,j+1}},\ \xi \ge  0,\ \xi \ \textrm{small enough}, \nonumber\\
& \ \bfm{h}_j := \left( \underbrace{0\dots,0}_{j-2},\const_{j,1},\const_{j,2},\const_{j,3},\underbrace{0\dots,0}_{2n-j-4}\right)^T\!\! .
\end{align}
%
Suppose now additionally ${j<n-1},\ {\bfm{t}^0 \in  \varisp{\fFunj{j+1}}.}$
By \lemnref{lem:fcatbSj} this is possible only if $\bfm{t}^0 \in {\smfp{j-1,j,j+1,j+2}}$, and we may use \eqref{eq:varabass} again, with $j \to j+1$.
An application of the Cauchy–Schwarz inequality
\begin{align}\label{eq:schjhjp1}
0 < \bfm{h}_j^T \bfm{h}_{j+1} & = \const_{j,2} \, \const_{j+1,1} + \const_{j,3} \, \const_{j+1,2} \le \\
&\le \sqrt{\const_{j,2}^2 + \const_{j,3}^2} \sqrt{ \const_{j+1,1}^2 + \const_{j,3}^2} < 
 \nrd{\bfm{h}_j} \nrd{\bfm{h}_{j+1}} \nonumber
\end{align}
verifies that directions of $\bfm{h}_j$ and $\bfm{h}_{j+1}$ are separated by a constant angle ${\psi_j \in \left(0,\frac{\pi}{2}\right)}$ since the inequality \eqref{eq:schjhjp1} is strict. 
Based upon this discussion we are able now to prove the major observation of this section.
\begin{lem}\label{lem:final}
Suppose the assumptions of \thmnref{theorem1} are satisfied. The functions $\fa, \fb$ can be modified to  $\tfa, \tfb$ in such a way that 
\begin{equation}\label{eq:fatilderel}
\sign\ \tfap{\vv{t}} = (-1)^n, \  \vv{t} \in \smfp{j-1}, \quad \sign \tfap{\qpp{j-1}} = (-1)^{n+1}, 
\end{equation}
and
\begin{equation}\label{eqfamoduni}
\sign\ \tfbp{\vv{t}} = 1, \ \vv{t} \in \smfp{n+j-1}, \quad \sign\ \tfbp{\qpp{n+j-1}} = -1.
\end{equation}
The substitution  ${\fa \to \tfa, \, \fb \to \tfb}$ doesn't bring about any additional solution of the system \eqref{system-pol-form}. 
\end{lem}
\begin{pf}
Let $0 < \charfunPara < \charfunParb \ll 1$, and let
\begin{equation}\label{eq:hidef}
\charfuni{\mathcal{W}}: \RR^{2n-2} \to \RR: \bfm{t} \mapsto \charfunNopp{\mathcal{W}}{\bfm{t}} := \charfunpp{\mathcal{W}}{\bfm{t}} := 
\left\{
\begin{matrix}
1, & \dist{ \bfm{t}}{\mathcal{W}} \le \charfunPara,\\
0, & \dist{ \bfm{t}}{\mathcal{W}} \ge \charfunParb,\\
0< \pika < 1, & \textrm{otherwise},
\end{matrix}
\right.
\end{equation}
denote a smooth characteristic wrapper function of a set $\mathcal{W} \subset\RR^{2n-2}$ with $\textrm{dist}$  denoting Hausdorff distance based upon Euclidean norm. Let us consider 
\eqref{eq:fatilderel} first.
\begin{figure}[htb]
	\centering
 \begin{minipage}{\textwidth}
    \centering \includegraphics[width=0.60\textwidth]{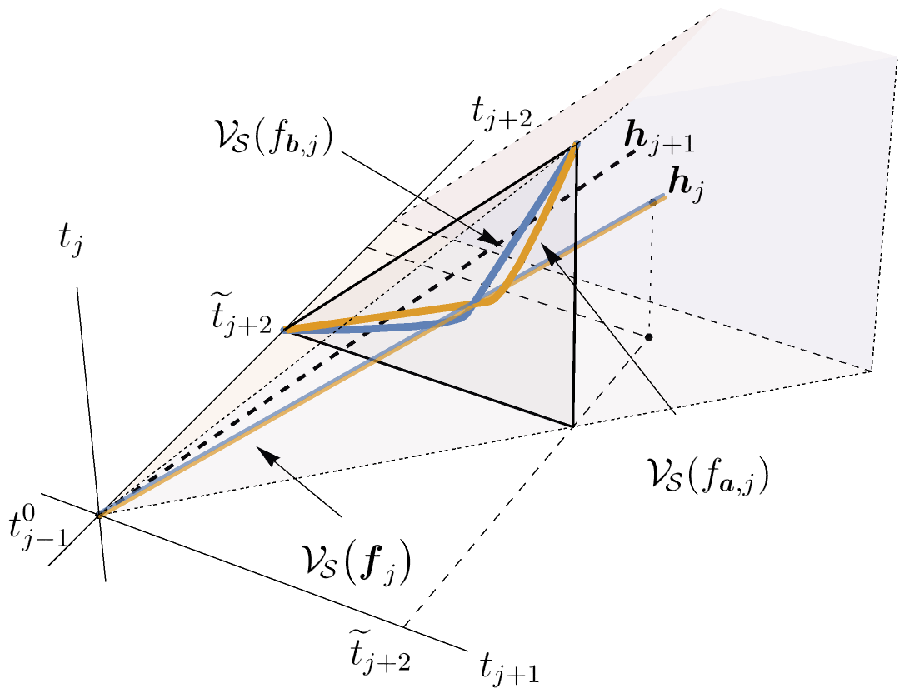}
 \end{minipage}
	\caption{The case $\Mi$ of $\smfp{j-1,j,j+1,j+2}$ in the proof of \lemnref{lem:fcatbSj}. The parameter  $\tti{j+2}$ is chosen small enough so that the varieties $\varisp{\fa}$ and $\varisp{\fb}$ could be expressed as functions, and $\varisp{\fFunj{j}}$ as a curve parameterized by $t_{j+2}$.}
\label{figSjmOnejjPTwo}
\end{figure}
From  \lemnref{lem:lemma1} and \lemnref{lem:fcatbSj} \ennref{en:fcatbSji} we conclude that there are only three boundary parts of 
$\smfp{j-1}$ that require a particular attention: 
\begin{enumerate}[label=(\roman*)]\setlength\itemsep{-4pt}
\item $\Mi =\smfp{j-1,j,j+1,j+2}$, 
\item $\Mii = \smfp{j-1,j,j+1} \setminus \Mi,$ 
\item $\Miii = \interior{\smfp{j-1,j-1+k}}$, with $k$ odd. 
\end{enumerate} 
If $\vtZ \in \smfp{j-1,j,j+1}$, \lemnref{lem:fcatbSj} implies that there exists 
$$
{  \tti{j+2} = \tti{j+2}\left(\tiZ{j-1}\right) ,\, 0 < \tti{j+2} < \coepsZeroiS{1},}
$$ 
such that the varieties $\varisp{\fa}$ and $\varisp{\fb}$ restricted to a tetrahedron $\triang$,
$$
\triang: \ \tiZ{j-1} \le \ti{j} \le \ti{j+1} \le \ti{j+2} \le \tti{j+2} \le \tiZ{j+3}, 
$$
are functions of variables $t_{j+1}$ and $t_{j+2}$ \figref{figSjmOnejjPTwo}. By \eqref{eq:varabass}, we may assume that $\varisp{\fFunj{j}}$
is a curve parameterized by $\ti{j+2} \in [\tiZ{j-1},\tti{j+2}]$. If not, a smaller $\tti{j+2} > 0$  can be found that allows this assumption.
\begin{figure}[htb]
	\centering
 \begin{minipage}{\textwidth}
    \centering \includegraphics[width=0.60\textwidth]{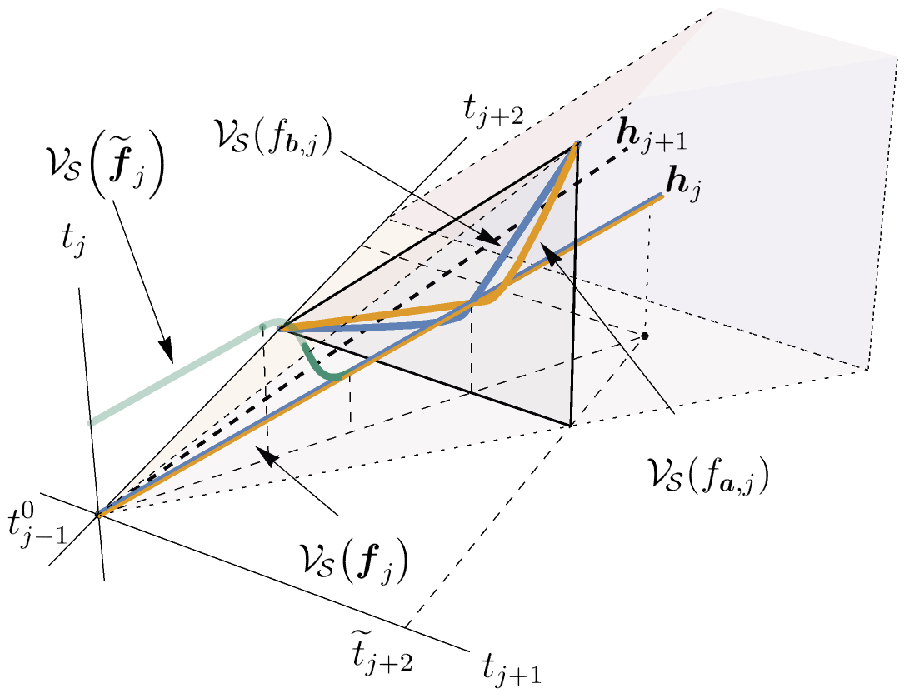}
 \end{minipage}
	\caption{The case $\Mi$ of $\subset\smfp{j-1,j,j+1,j+2}$ in the proof of \lemnref{lem:fcatbSj}. The function $\fa$ is modified to$\tfa$. As a consequence, $\varisp{\fFunj{j}}$ is modified to $\varisp{\tfFunj{j}}$ but no additional intersection with $\varisp{\fFunj{j+1}}$ is introduced.}
\label{figSjmOnejjPOne}
\end{figure}
Let us now choose somewhat arbitrary
\begin{align*}
\bfm{x} := \left(x_i\right)^T := - \vt^0 + \varisp{\fFunj{j}}_{\big|_{\ti{j+2}=\tti{j+2}}},\quad \charfunParb = \frac{2}{3} \nrd{\bfm{x}}, \quad \charfunPara = \frac{2}{3}\charfunParb, \ \varrho := x_j.
\end{align*}
A modification
$$
\fa \to \tfa := \fa + \left(-1\right)^n \varrho \, \charfuni{\smfp{j-1,j,j+1}}
$$
by \lemnref{lemmodstep} doesn't introduce  any additional solution of \eqref{system-pol-form} provided the corresponding modified part of
$\varisp{\tfFunj{j}}$ doesn't  meet $\varisp{\fFunj{j+1}}$. This is obvious for the set $\Mii$, and also by \eqref{eq:varabass} and \eqref{eq:schjhjp1} for the case $\Mi$ if $\tti{j+2}$ is small enough \figref{figSjmOnejjPOne}. However, $\smfp{j-1,j,j+1}$ is compact, so the smallest values $\tti{j+2}, \charfunParb$ and $\varrho$ can be found, and the modification \eqref{eqfamoduni} works uniformly. Note that this modification covers also the points of $\Miii$ that are close enough to 
$\smfp{j-1,j,j+1}$. For those that are not, consider ${\vt^0 \in \varisp{\tfa} \cup \boundary{\smfp{j-1} }}$. \lemnref{lem:fcatbSj} implies that in small enough neighbourhood 
$$
\left\{ \vt \in \interior{\sm} \, \big| \, \nrd{\vt^0-\vt}{} < \varepsilon\left(\vt^0\right) \right\}
$$
the varieties $\varisp{\tfa}$ and $\varisp{\fb}$ are separated. Since $\varisp{\tfa} \cup \boundary{\smfp{j-1} }$ is compact, the smallest $\varepsilon$ exists. So one may again construct a modification of the type \eqref{eq:hidef} that will by \lemnref{lemmodstep} introduce no additional solutions of \eqref{system-pol-form}. By combining both wrappers we finally obtain the admissible modification of $\fa$ at $\smfp{j-1}$
\begin{equation}\label{eq:bothwrap}
\fa \to \tfa := \fa + \left(-1\right)^n \varrho^*_{\bfm{a},j} \, \charfuni{\smfp{j-1}}, \quad \varrho^*_{\bfm{a},j} > 0,
\end{equation}
that removes also $\vara \backslash \varisp{\fa}$ from the simplex boundary and proves the first part of \eqref{eq:fatilderel}. The first part of \eqref{eqfamoduni} follows in a similar way. We have to verify only that a modification $\fa \to \tfa$ doesn't spoil the arguments used in the proof of \eqref{eq:fatilderel}. At the common face of $\smfp{j-1}$ and $\smfp{n+j-1}$ we observe
$$
\smfp{j-1,n+j-1} \cup \varisp{\fc} = \emptyset
$$
if $n$ is even since 
$$
\sign \fc\left(\vt\right)   = (-1)^{n} = 1 = \sign \fc(\,\widetilde{\vt}\,), \quad \vt \in \smfp{j-1}, \widetilde{\vt} \in \smfp{n+j-1}.
$$
If $n$ is odd, modifications $\tfa - \fa, \tfb-\fb $ are of opposite sign, so also the first assertion of \eqref{eqfamoduni} is confirmed.
Suppose now that $ \fa \to \tfa$, $\fb \to \tfb$ has been carried out as in \eqref{eq:bothwrap}, so $\sign\tfbp{\vt} = 1, \ \vt \in \smfp{n+j-1}$. So one may $\tfa$ by \lemnref{lem:fcatbSj} further modify
\begin{equation}\label{eq:modagain}
\tfa \to \tfa + \left(-1\right)^{n+1} \varrho^*_{\bfm{a},n+j-1} \, \charfuni{\qpp{n+j-1}},
\end{equation}
for sufficiently small $\varrho^*_{\bfm{a},n+j-1} > 0$. The modification \eqref{eq:modagain} doesn't introduce any new intersection between $\tfa$, and $\tfb$. The argument works at $\qpp{n+j-1}$ too. This completes the proof.
\qed
\end{pf}

%
\section{Proof of Theorem~\ref{theorem1}}\label{sec:Proof of Theorem1}

Let us apply \citep[][Theorem 2.1]{Vrahatis-2016} on 
$$
\tfFun := \left( \tfunpp{\bfm{a}}{1}, \tfunpp{\bfm{b}}{1}, \dots,  \tfunpp{\bfm{a}}{n-1}, \tfunpp{\bfm{b}}{n-1} \right)^T,
$$
derived from $\fFun$  by \lemnref{lem:final}.
The choice of the map $\lstar$ defined in \eqref{lstar} that determines the simplex-vertex pairs at which the functions $\tfc$ are of different sign is straightforward,
$$
 \lstarp{{\vv{a},j}} = j-1, \quad  \lstarp{{\vv{b},j}} = n + j - 1,\quad j = 1,2,\dots,n-1.
$$ 
By \lemnref{lem:final}, 
\begin{align*}
& \sign \tfap{\vt}  = (-1)^n \neq \sign \tfap{\qpp{j-1}} = (-1)^{n+1}, \quad \vt \in \smfp{j-1}, \nonumber \\
& \sign \tfbp{\vt} = 1  \neq \sign \tfbp{\qpp{n+j-1}} = -1, \quad \vt \in \smfp{n+j-1}, \nonumber
\end{align*}
and \eqref{pog1} is satisfied. With $\ldstar = n - 1$, it remains to verify \eqref{pog2},
\begin{equation}\label{eq:pogzvver}
 \sign \tfFunp{\vt}  \neq \sign \tfFunp{\qpp{n-1}}, \vt \in \smfp{n-1}.
\end{equation}
Since $1\le j < n$, and $\qpp{n-1} \in \smfp{j-1,n+j-1}$, we have
$$
\sign\ \tfap{\qpp{n-1}} = (-1)^n, \ \sign\ \tfbp{\qpp{n-1}} = 1,\quad j= 1,\dots,n-1.
$$
\lemnref{lem:lemma1} further reveals
$$
\sign\ \fcp{\vt} = (-1)^j,\quad \vt \in \interior{\smfp{n-1}}.
$$
For an odd $n-j$ then
$$
\sign\ \tfap{\vt} \ne \sign\ \tfap{\qpp{n-1}},\quad \vt \in \smfp{n-1} \backslash \supp\charfuni{\smfp{j-1}}.
$$
Note that points in $\smfp{n-1} \cap \supp\charfuni{\smfp{j-1}}$ are arbitrary close to $\smfp{j-1,n-1}$. Similarly,
for an odd $j$, we conclude
$$
\sign\ \tfbp{\vt} \ne \sign\ \tfbp{\qpp{n-1}},\quad \vt \in \smfp{n-1} \backslash \supp\charfuni{\smfp{n+j-1}}.
$$
\begin{table}[htb]
\begin{center}
\begin{tabular}{ |c|c|c|c|c|c|c|}
\hline
 &  &  &  $\sign\fa$ at   & $\sign\tfa$ at & $\sign\fb$ at   & $\sign\tfb$ at  \\
$n$ & $j$ & $n-j$ &   $\interior{\smfp{n-1}}$ & $\qpp{n-1}$ &    $\interior{\smfp{n-1}}$ & $\qpp{n-1}$ \\
\hline
 odd & odd  & even & $-1$  &  $-1$ &  $-1$ & $1$ \\
 odd & even & odd  &  $\ \ \  1$  &  $-1$ &  $\ \ \ 1$ &  $1$\\
 \hline
 even & odd & odd  &  $-1$ &  $\  \ \ 1$ &  $-1$ &  $1$\\
even& even & even &  $\ \ \ 1$  &  $\ \ \  1$ &  $\  \ \ 1$ &  $1$\\
\hline
\end{tabular}
\caption {Sign distribution of  $\fa$, $\tfa$,  $\fb$, $\tfb$ at $\smfp{n-1}$ and $\qpp{n-1}$. }
\label{tab:signdis}
\end{center}
\end{table}
Tab.~\ref{tab:signdis} summarizes the sign distribution detected. If $n$ is odd, the pairs $\tfa$ and $\tfb$ demonstrate similar behavior  regardless of $j$ been even or odd. Suppose it is odd. Then there is no $\varisp{\tfb} \cap \smfp{n-1}$ except near ${\smfp{n+j-1} \cap \smfp{n-1}}$.
On the other hand, $\sign\tfap{\qpp{j-1}} = 1$, and $\sign\tfap{\qpp{j-1}} \in \{0,1\}$ except near ${\smfp{j-1} \cap\smfp{n+j-1}}$. So 
$\tfFunj{j}$ confirms \eqref{eq:pogzvver} except possibly near 
 $$
 \smfp{j-1} \cap \smfp{n-1} \cap \smfp{n+j-1} = \smfp{j-1,n-1,n+j-1}.
 $$
 But 
 $$
\bigcap\limits_{j=1}^{n-1} \smfp{j-1,n-1,n+j-1} = \emptyset,
 $$
 what confirms \eqref{eq:pogzvver} for an odd $n$. If $n$ is even, $\tfFunj{j}$  confirms again \eqref{eq:pogzvver} except possibly near 
 $$
 \smfp{j-1} \cap \smfp{n-1} \cap \smfp{n+j-1} = \smfp{j-1,n-1,n+j-1}.
 $$
 To observe this note that $\tfa$ verifies \eqref{eq:pogzvver} except close to $\smfp{j-1} \cap \smfp{n-1}$ for an odd $j$, and close to $\smfp{n+j-1} \cap \smfp{n-1}$ for an even one. The rest of the proof is obvious. This finally completes the proof of  Theorem~\ref{theorem1}. 
\hfill\qed\break

A brief inspection of the proofs of Theorem~\ref{theorem1} and the lemmas of Section~\ref{sec:Analysis of polynomials}   it is based upon reveals that there is no need for the constants ${\deap{\ell}, \debp{\ell}}$ to be the same for distinct values of $j$ as long as their required sign behavior is preserved. The following corollary  illuminates this fact.
\begin{cor}\label{cor:corollary2} 
Suppose that the data \eqref{remdat} that determine $\fa$ and $\fb$ depend on $j$, 
$$
\deap{\ell} =\deap{\ell,j}, \ \debp{\ell} = \debp{\ell,j}, \quad \ell = j-1,j,\dots,n+j-1.
$$
Suppose that the following assumptions are fulfilled for each $j, j = 1,2,\dots,n-1$, separately:
the data differences
\begin{gather*}
\deap{\ell,j} , \ \debp{\ell,j}, \quad  \ell = j-1,j,\dots,n+j-1, 
\end{gather*}
are all positive or all negative, and   determinants
\begin{gather*}
\dt{\deap{\ell,j}}{\deap{\ell+1,j}}{\debp{\ell,j}}{\debp{\ell+1,j}} ,
\quad \ell = j-1,j, \dots,n+j-2, 
\end{gather*}
are all of the same sign. The system of equations \eqref{system-pol-form} has at least one  solution.
\end{cor}

\section{Proof of Theorem~\ref{theorem2}}\label{sec:Proof of Theorem 2}

We will prove Theorem~\ref{theorem2} with help of Remark~\ref{rem:remark1} and Corollary~\ref{cor:corollary2}. Without losing generality we may assume that the sign of determinants \eqref{eqthmsupa} is positive. Consider fixed ${j, 1\le j \le n -1}$. By  Remark~\ref{rem:remark1}, we may modify the data  \eqref{remdat} by a matrix product
$$
M \begin{pmatrix} 
\debp{n+j-1} & - \deap{n+j-1}  \\
- \debp{j-1}  &\deap{j-1}
\end{pmatrix},
\quad M :=
 \begin{pmatrix} 
1 & \varepsilon \\
\varepsilon  & 1
\end{pmatrix}, \quad \varepsilon  > 0.
$$
This yields the new data differences
\begin{equation*}
\begin{pmatrix} 
\detap{\ell}  \\
\detbp{\ell} 
\end{pmatrix}
:=
M
\begin{pmatrix} 
\debp{n+j-1} & - \deap{n+j-1}  \\
- \debp{j-1}  &\deap{j-1}  
\end{pmatrix}
\begin{pmatrix} 
\deap{\ell}  \\
\debp{\ell} 
\end{pmatrix},
\quad \ell = j-1,j,\dots,n+j-1.
\end{equation*}
But
\begin{align*}
\detap{\ell} & = \dt{\deap{\ell}}{\deap{n+j-1}}{\debp{\ell}}{\debp{n+j-1}} + \varepsilon  \dt{\deap{j-1}}{\deap{\ell}}{\debp{j-1}}{\debp{\ell}} > 0, \nonumber \\
&\quad\quad\quad  \quad \quad\quad\quad  \quad  \quad\quad\quad  \quad \quad\quad\quad  \quad 
\quad\quad
\ell = j-1,j,\dots,n+j-1,\\
\detbp{\ell} & =  \varepsilon \dt{\deap{\ell}}{\deap{n+j-1}}{\debp{\ell}}{\debp{n+j-1}} +  \dt{\deap{j-1}}{\deap{\ell}}{\debp{j-1}}{\debp{\ell}} > 0, \nonumber \\
\end{align*}
by the assumption of  Theorem~\ref{theorem2}.
Also, 
\begin{align*}
\dt{\detap{\ell}}{\detap{\ell+1}}{\detbp{\ell}}{\detbp{\ell+1}} = & (1-\varepsilon^2)  \dt{\deap{j-1}}{\deap{n+j-1}}{\debp{j-1}}{\debp{n+j-1}} \dt{\deap{\ell}}{\deap{\ell+1}}{\debp{\ell}}{\debp{\ell+1}} > 0, \\
& \quad\quad\quad\quad\quad\quad\quad\quad \quad\quad \ell = j-1,j,\dots,n+j-2,
\end{align*}
for $\varepsilon$ small enough. But then the assumptions of Corollary~\ref{cor:corollary2} are met. This concludes the proof of the theorem. \qed

\clearpage


\bibliographystyle{elsarticle-num} 
%

\begin{thebibliography}{10}
\expandafter\ifx\csname url\endcsname\relax
  \def\url#1{\texttt{#1}}\fi
\expandafter\ifx\csname urlprefix\endcsname\relax\def\urlprefix{URL }\fi
\expandafter\ifx\csname href\endcsname\relax
  \def\href#1#2{#2} \def\path#1{#1}\fi

\bibitem{deBoor-Hoellig-Sabin-87-High-Accuracy}
C.~de~Boor, K.~H{\"o}llig, M.~Sabin, High accuracy geometric {H}ermite
  interpolation, Comput. Aided Geom. Design 4~(4) (1987) 269--278.

\bibitem{HolligKoch-1996}
K.~H{\"o}llig, J.~Koch, Geometric {H}ermite interpolation with maximal orderand
  smoothness, Comput. Aided Geom. Design 13~(8) (1996) 681--695.

\bibitem{HolligKoch-1995}
K.~H{\"o}llig, J.~Koch, Geometric {H}ermite interpolation, Comput. Aided Geom.
  Design 12~(6) (1995) 567--580.

\bibitem{MorkenScherer-1997}
K.~M{\o}rken, K.~Scherer, A general framework for high-accuracy parametric
  interpolation, Math. Comp. 66~(217) (1997) 237--260.

\bibitem{Scherer-2000}
K.~Scherer, Parametric polynomial curves of local approximation order {$8$},
  in: Curve and Surface Fitting (Saint Malo, 1999), Vanderbilt Univ. Press,
  Nashville, TN, 2000, pp. 375--384.

\bibitem{KozakZagar-2004}
J.~Kozak, E.~{\v{Z}}agar, On geometric interpolation by polynomial curves, SIAM
  J. Numer. Anal. 42~(3) (2004) 953--967.

\bibitem{JaklicKozakKrajncZagar-2007b}
G.~Jakli\v{c}, J.~Kozak, M.~Krajnc, E.~{\v{Z}}agar, On geometric interpolation
  by planar parametric polynomial curves, Math. Comput 76~(260) (2007)
  1981--1993.

\bibitem{Brysiewicz-2021}
T.~Brysiewicz, Necklaces count polynomial parametric osculants, Journal of
  Symbolic Computation 103 (2021) 95--107.

\bibitem{Vrahatis-2016}
M.~N. Vrahatis, Generalization of the {B}olzano theorem for simplices, Topology
  and its Appl. 202 (2016) 40--46.

\bibitem{LachanceSchwartz-1991}
M.~Lachance, A.~Schwartz, Four point parabolic interpolation, Comput. Aided
  Geom. Design 8~(2) (1991) 143--150.

\bibitem{Morken-1995}
K.~M{\o}rken, Parametric interpolation by quadratic polynomials in the plane,
  in: M.~D{\ae}hlen, T.~Lyche, L.~L. Schumaker (Eds.), Mathematical Methods for
  Curves and Surfaces, Vanderbilt University Press, Nashville, USA, 1995, pp.
  385--402.

\bibitem{KozakKrajnc-2007}
J.~Kozak, M.~Krajnc, Geometric interpolation by planar cubic polynomial curves,
  Comput. Aided Geom. Design 24~(2) (2007) 67--78.

\bibitem{Krajnc-2009}
M.~Krajnc, Geometric {H}ermite interpolation by cubic {$G^1$} splines,
  Nonlinear Anal. 70~(7) (2009) 2614--2626.

\bibitem{Degen-approximation-95}
W.~L.~F. Degen, High accuracy approximation of parametric curves, in:
  Mathematical methods for curves and surfaces ({U}lvik, 1994), Vanderbilt
  Univ. Press, Nashville, TN, 1995, pp. 83--98.

\bibitem{deBoor-Guide-2001}
C.~de~Boor, A practical guide to splines, revised Edition, Vol.~27 of Applied
  Mathematical Sciences, Springer-Verlag, New York, 2001.

\end{thebibliography}
%
%


\end{document}